\documentclass[letterpaper, 12pt]{amsart}

\usepackage[all]{xy}
\usepackage{latexsym}
\usepackage{amssymb}
\usepackage{fullpage}
\usepackage{amsmath}
\usepackage{amsthm}
\usepackage[dvips]{graphicx}
\usepackage{yfonts}

\newtheorem{theorem}{Theorem}[section]
\newtheorem{lemma}[theorem]{Lemma}
\newtheorem{corollary}[theorem]{Corollary}
\newtheorem{defn}[theorem]{Definition}
\newtheorem{proposition}[theorem]{Proposition}

\newtheorem{question}[theorem]{Question}
\newtheorem{remark}[theorem]{Remark}
\numberwithin{equation}{section}

\newcommand{\V}{\Vert}

\newcommand{\la}{\langle}
\newcommand{\ra}{\rangle}

\title{Perturbations of completely positive maps and strong NF algebras}
\author{Caleb Eckhardt}
\address{Department of Mathematics, University of Illinois, Urbana, IL, 61801}
\email[Caleb Eckhardt]{ceckhard@uiuc.edu}
\date{}

\begin{document}
\maketitle
\begin{abstract}
Let $\phi:M_n\rightarrow B(H)$ be an injective, completely positive contraction with $\V\phi^{-1}:\phi(M_n)\rightarrow M_n\V_{cb}\leq1+\delta(\epsilon).$  We show that if either (i) $\phi(M_n)$ is faithful modulo the compact operators or (ii) $\phi(M_n)$ approximately contains a rank 1 projection, then there is a complete order embedding $\psi:M_n\rightarrow B(H)$ with $\V\phi-\psi\V_{cb}<\epsilon.$   We also give examples showing that such a perturbation does not exist in general.  As an application, we show that every $C^*$-algebra $A$ with $\mathcal{OL}_\infty(A)=1$ and a finite separating family of primitive ideals is a strong NF algebra, providing a partial answer to a question of Junge, Ozawa and Ruan. 
\end{abstract}
\section{Introduction}
The main goal of this paper is to provide a clearer picture of the relationship between Blackadar and Kirchberg's strong NF algebras and Junge, Ozawa and Ruan's notion of $\mathcal{OL}_\infty$ structure. Our main technical tools in this endeavor are perturbations of completely positive maps.  We feel that these perturbation results may have applications in other areas of operator algebras, so we first describe them in some detail.

Let $B(H)$ denote the space of bounded linear operators on a Hilbert space $H$, and let $M_n$ denote $B(H)$ when the dimension of $H$ is $n.$ Following \cite[4.1]{Blackadar97} we say  that a linear map $\phi$ between $C^*$-algebras is a \emph{complete order embedding} if $\phi$ is a completely positive, complete isometry.

We start with the following question:
\begin{question} \label{ques:Chapter2question} Let $\epsilon>0$ be given. Does there exist a $\delta>0$ such that if $\phi:M_n\rightarrow B(H)$ is an injective, completely positive contraction  with 
$\V \phi^{-1}:\phi(M_n)\rightarrow M_n\V_{cb}<1+\delta$, then there is a complete order embedding $\psi:M_n\rightarrow B(H)$ such that $\V\phi-\psi\V<\epsilon$? In particular, can $\delta$ be chosen independent of $n$?
\end{question}
David Kerr and Hanfeng Li showed in \cite[Lemma 7.1]{Kerr08} that if $\epsilon$ depends on both  $\delta$ and $n$, then one can  
answer Question \ref{ques:Chapter2question} affirmatively. We show (Section \ref{sec:dimensiondep}) that the dependence on $n$ in Kerr and Li's lemma cannot be omitted.  Hence the answer to Question \ref{ques:Chapter2question} is no, in general.

 On the other hand we show  that in two extreme cases, the answer to Question \ref{ques:Chapter2question} is yes. First we show that if the range of $\phi$ is faithful modulo $K(H)$ (the compact operators on $H$), then we can perturb independent of $n$ (Theorem \ref{thm:antiliminalperturb}).  At the other end of the spectrum, we show that if $\phi$ approximately contains a rank 1 projection, then we can  perturb $\phi$ independent of $n$ (Theorem \ref{thm:minimalperturb}).

These two special cases prove to be the most applicable, due to the well-known fact that for any irreducible representation $(\pi, H)$ of a $C^*$-algebra $A$ we either have $\pi(A)\cap K(H)=\{0\}$ or $K(H)\subseteq \pi(A).$

In the second part of the paper, we apply our perturbation results to the study of strong NF algebras. Blackadar and Kirchberg introduced and studied strong NF algebras in \cite{Blackadar97, Blackadar01, Blackadar07}. A $C^*$-algebra $A$ is called strong NF if for every finite subset $F\subset A$ and $\epsilon>0$, there is a finite dimensional $C^*$-algebra $B$ and a  complete order embedding  $\psi:B\rightarrow A$ such that $\max_{x\in F}dist(x,\psi(B))<\epsilon.$
Strong NF algebras form a large and natural class of $C^*$-algebras. Indeed, it was shown in \cite{Blackadar01} that strong NF algebras are precisely the nuclear, inner quasidiagonal $C^*$-algebras. In particular, every simple, nuclear, quasidiagonal $C^*$-algebra is strong NF.

In \cite{Junge03}, the authors generalized the notion of strong NF algebras by  introducing the invariant $\mathcal{OL}_\infty(\cdot)$ for  $C^*$-algebras.  Let $A$ be a $C^*$-algebra and $\lambda>1.$  $A$ is called an $\mathcal{OL}_{\infty, \lambda}$ space if 
  for every finite subset $F\subset A$, there is a finite dimensional $C^*$-algebra $B$ and an injective linear map $\phi:B\rightarrow A$ such that $F\subset \phi(B)$ and $\V\phi\V_{cb}\V\phi^{-1}:\phi(B)\rightarrow B\V_{cb}<\lambda.$ 
One then defines $\mathcal{OL}_\infty(A)$ to be the infimum of all $\lambda$ such that $A$ is an $\mathcal{OL}_{\infty,\lambda}$ space. It was shown in \cite{Junge03} that $\mathcal{OL}_\infty(A)$ is finite precisely when $A$ is nuclear. 

In this work we are interested in $C^*$-algebras $A$ with $\mathcal{OL}_\infty(A)=1.$  It was shown in \cite{Junge03} that all $C^*$-algebras with $\mathcal{OL}_\infty(A)=1$ are (nuclear) and quasidiagonal.  On the other hand in \cite{Eckhardt08}, examples were given   of nuclear quasidiagonal $C^*$-algebras $A$ with $\mathcal{OL}_\infty(A)>1.$ At present we do not have a good description  of those $C^*$-algebras with $\mathcal{OL}_\infty(A)=1.$  A goal of this paper is to bring us closer to a ``good'' description of these algebras.

A simple perturbation argument shows that all strong NF algebras $A$ have $\mathcal{OL}_\infty(A)=1,$ and it was asked in \cite{Junge03} whether or not all $C^*$-algebras with $\mathcal{OL}_\infty(A)=1$  are strong NF algebras. 

 We use the perturbation arguments of the first part of the paper to show that if $A$ has a finite, separating family of primitive ideals and $\mathcal{OL}_\infty(A)=1,$ then $A$ is a strong NF algebra (Theorem \ref{thm:maintheoremin3}). The general case is still open, and we discuss the  known necessary conditions at the end of the paper. We then apply Theorem \ref{thm:maintheoremin3} to show that $\mathcal{OL}_\infty$ is not continuous with respect to inductive limits and not multiplicative with respect to tensor products, even when one of the algebras is AF.

We use the following shorthand notation throughout the paper: Let $\phi:A\rightarrow B$ be a linear map.  If $\phi$ is a unital completely positive map, we say that $\phi$ is a UCP map.  If $\phi$ is a completely positive contraction, we say that $\phi$ is a CPC. We define
\begin{equation*}
\phi^{(n)}:=id_{M_n}\otimes \phi: M_n\otimes A\rightarrow M_n\otimes B.
\end{equation*} 
If $\phi$ is injective with closed range, we define
\begin{equation*}
 \V \phi^{-1}\V_{cb}:=\V \phi^{-1}:\phi(A)\rightarrow A\V_{cb}.
\end{equation*}
We write $\ell^2(n)$ for $n$-dimensional Hilbert space with standard basis vectors $e_1,...,e_n$ and $(e_{ij})_{i,j=1}^n$ as standard matrix units for $M_n.$ We let $A^+$ denote the positive cone of the $C^*$-algebra $A$ and for $x,y\in A$ we write $[x,y]=xy-yx.$

\section{Perturbations} 
The objects of study in this section are injective CPCs, $\phi:M_n\rightarrow B(H)$ with $\V\phi^{-1}\V_{cb}$ close to 1.  We first answer Question \ref{ques:Chapter2question} affirmatively when dim$(H)=n.$ Then we gather some non-unital analogs of well-known  theorems about UCP maps.  We finish this section with the general answer to Question \ref{ques:Chapter2question} discussed in the introduction.

\subsection{The case where dim$(H)=n$}
\begin{lemma} \label{lem:perturbisometries} Let $0<\epsilon<1/2$, $H, K$ be Hilbert spaces and $v:H\rightarrow K$ an isometry.  Suppose that there is a projection
 $p\in B(K)$ such that $\V vv^*-p\V<\epsilon.$  Then there is an isometry $w:H\rightarrow K$ such that
$\V v-w\V<2\epsilon$ and $ww^*=p.$
\end{lemma}
\begin{proof} Consider the polar decomposition of $pv=w|pv|.$  It follows from basic spectral theory that $w$ has the desired properties.
\end{proof}
 The following proposition (in the unital case) can be deduced from Eric Christensen's paper \cite{Christensen77}.  It appears that proving it directly versus applying \cite[Lemma 3.3]{Christensen77} requires the same amount of work, so we prove it directly.

\begin{proposition} \label{prop:samerange} Let $57^{-1/2}>\delta>0$ and
$\phi:M_n\rightarrow M_n$ a CPC with $\V \phi^{-1}\V<1+\delta$. Then there is a *-automorphism $\pi$ of $M_n$ such that
$\V\pi-\phi\V_{cb}<57\sqrt{\delta}.$
\end{proposition} 
\begin{remark} We emphasize that the condition $\V \phi^{-1}\V<1+\delta$ in Proposition \ref{prop:samerange} is not a typo, i.e. it is not necessary that the cb-norm of the inverse is  close to 1, only that the inverse of the 1 norm is close to 1. This won't be particularly important for us, but it does have the happy consequence of less notation in Lemma \textup{\ref{lem:cruxlemma}}.
\end{remark}

\begin{proof} Suppose first that $\phi(1)=1.$ Let $(\sigma, H, v:\ell^2(n)\rightarrow H)$ be the Stinespring dilation of $\phi.$  
Since $\phi(1)=1,$  $v$ is an isometry.

We first show that the $C^*$-algebra $\sigma(M_n)$ approximately commutes with the projection $vv^*.$ To this end, let $u\in M_n$ be unitary.
Then $\V \phi^{-1}(u)\V\leq 1+\delta.$  Hence,
\begin{align}
 (1+\delta)^2 & \geq  \V vv^*\sigma(\phi^{-1}(u))\V^2 \notag \\
&= \V vv^*\sigma(\phi^{-1}(u))vv^*+vv^*\sigma(\phi^{-1}(u))(1-vv^*)\V^2 \notag \\
&= \V vuv^*+vv^*\sigma(\phi^{-1}(u))(1-vv^*)\V^2 \notag \\
&= \V vv^*+vv^*(\sigma(\phi^{-1}(u))(1-vv^*)\sigma(\phi^{-1}(u))^*vv^*\V  \notag \\
&=1+\V vv^*(\sigma(\phi^{-1}(u))(1-vv^*)\V^2, \label{align:cornercomp}
\end{align}
where the last line follows by spectral theory performed in the $C^*$-algebra $vv^* B(H) vv^*.$
By (\ref{align:cornercomp}) applied to both $u$ and $u^*$  we obtain,
\begin{equation}
 \V [ \sigma(\phi^{-1}(u)),vv^*]\V \leq (2\delta+\delta^2)^{1/2}\leq 2\sqrt{\delta}. \label{eq:approxcommute}
\end{equation}
Applying the Russo-Dye theorem to $\phi(u)$ we obtain unitaries $v_1,...,v_r\in M_n$ and positive scalars $\lambda_1,...,\lambda_r$ that sum to 1 such that
\begin{equation*}
 u=\sum_{i=1}^r \lambda_i \phi^{-1}(v_i).
\end{equation*}
By (\ref{eq:approxcommute}) applied to $v_1,...,v_r$ it follows that 
\begin{equation*}
 \V [\sigma (u), vv^*]\V \leq 2\sqrt{\delta}\quad \textrm{ for every unitary }u\in M_n.
\end{equation*}

Let $\mu$ denote normalized Haar measure on $\mathcal{U}_n$, the unitary group of $M_n.$ Define
\begin{equation*}
x=\int_{\mathcal{U}_n}\sigma(u)(vv^*)\sigma(u)^* d\mu(u).
\end{equation*}
Then, $x\in\sigma(M_n)'$ and 
\begin{equation*}
\V x-vv^*\V\leq \int_{\mathcal{U}_n} \V\sigma(u)(vv^*)\sigma(u)^*-vv^*\V d\mu(u)\leq 2\sqrt{\delta}.
\end{equation*}
Now, $x$ may not be a projection, but it is close to one.  It is clearly positive, and
\begin{align}
\V x^2-x\V & \leq  \V x^2-(vv^*)\V+\V x-vv^*\V \notag \\
&\leq \V x^2-xvv^*\V +\V xvv^*-vv^*\V +\V x-vv^*\V \notag \\
&\leq 3(2\sqrt{\delta})=6\sqrt{\delta}. \label{align:speccalc}
\end{align}
Since $0\leq x\leq1,$ and by (\ref{align:speccalc}), it follows that
\begin{equation*}
\textrm{sp}(x)\subset [0,12\sqrt{\delta}]\cup [1-12\sqrt{\delta},1].
\end{equation*}
Let $p$ be the spectral projection of $x$ associated with $[1-12\sqrt{\delta},1].$
Then, $p\in\sigma(M_n)'$ and
$\V p-x\V\leq 12\sqrt{\delta},$ hence
\begin{equation*}
\V p-vv^*\V\leq 14\sqrt{\delta}.
\end{equation*}
We apply Lemma \ref{lem:perturbisometries} to obtain an isometry $w:\ell^2(n)\rightarrow H$ such that
$\V w-v\V\leq 28\sqrt{\delta}$ and $ww^*=p.$  Now consider the map
\begin{equation*}
\pi:M_n\rightarrow M_n\textrm{ defined by }\pi(x)=w^*\sigma(x)w.
\end{equation*}
Then, $\V \pi-\phi\V_{cb}\leq 56\sqrt{\delta},$ and since $ww^*$ commutes with $\sigma(M_n)$ it follows that $\pi$ is a *-homomorphism. 

For the non-unital case, notice that $\phi^{-1}(1)\phi^{-1}(1)^*\leq (1+\delta)^2.$  Then
\begin{equation*}
 1= \phi(\phi^{-1}(1))\phi(\phi^{-1}(1))^*\leq \phi(\phi^{-1}(1)\phi^{-1}(1))^*\leq (1+\delta)^2\phi(1).
\end{equation*}
By basic spectral theory, it follows that $\V \phi(1)^{-1/2}-1\V\leq\delta.$ Hence defining
\begin{equation*}
 \psi(x)=\phi(1)^{-1/2}\phi(x)\phi(1)^{-1/2}
\end{equation*}
It follows that $\psi$ is unital and $\V \psi-\phi\V_{cb}\leq 2\delta<\sqrt{\delta}.$  We apply the first part of the proof to $\psi$ to obtain the conclusion.

\end{proof}
Eventually, we will want to perturb maps where the dimension of $H$ is arbitrary.  Our method will be to first cut down by a rank $n$-projection and then perturb the cutdown map via Theorem \ref{prop:samerange}.  The next lemma provides the justification for this method.
\begin{defn}
 Let $\phi:A\rightarrow B$ be a linear map between $C^*$-algebras and let $p\in B$ be a projection.  We define the map $\phi_p:A\rightarrow B$ as
\begin{equation*}
 \phi_p(x)=p\phi(x)p.
\end{equation*}

\end{defn}

\begin{lemma} \label{lemma:cutdownalmostiso}
 Let $\phi:M_n\rightarrow B(H)$ be a CPC.  Suppose there is a rank $n$ projection $p\in B(H)$ such that
$\phi_p$ is injective with $\V \phi_p^{-1}\V\leq 1+\delta.$  Then there is a complete order embedding $\psi:M_n\rightarrow B(H)$ such that $\V\psi-\phi\V_{cb}\leq 68\delta^{1/4}.$ Moreover $\psi=\psi_p+\psi_{(1-p)}$ where $\psi_p$ is a nonzero *-homomorphism.
\end{lemma}

\begin{proof}
 By Proposition \ref{prop:samerange} there is a *-homomorphism $\sigma:M_n\rightarrow pB(H)p$ such that 
\begin{equation}
\V\sigma-\phi_p\V_{cb}\leq 57\sqrt{\delta}. \label{eq:pphippert}
\end{equation}
Let $k\in\mathbb{N}$ be arbitrary.  Let $u\in M_k\otimes M_n$ be a unitary.  Then
\begin{align*}
 1&\geq \V \phi_p^{(k)}(u)+(1_k\otimes p)\phi^{(k)}(u)(1_k\otimes 1_n-1_k\otimes p)\V\\
&\geq \V \sigma^{(k)}(u)+(1_k\otimes p)\phi^{(k)}(u)(1_k\otimes 1_n-1_k\otimes p)\V-57\sqrt{\delta}\\
&=\V 1_k\otimes p+(1_k\otimes p)\phi^{(k)}(u)(1_k\otimes 1_n-1_k\otimes p)\phi^{(k)}(u)^*(1_k\otimes p)\V^{1/2}-57\sqrt{\delta}\\
&=\Big(1+\V(1_k\otimes p)\phi^{(k)}(u)(1_k\otimes 1_n-1_k\otimes p)\phi^{(k)}(u)^*(1_k\otimes p)\V\Big)^{1/2}-57\sqrt{\delta},
\end{align*}
where the last line follows by spectral theory. After rearranging and then applying the Russo-Dye theorem, it follows that
\begin{equation}
 \V (1_k\otimes p)\phi^{(k)}(x)(1_k\otimes 1_n-1_k\otimes p)\V\leq \sqrt{115}\delta^{1/4}  \label{eq:cornercomp2}
\end{equation}
for every $x\in M_k\otimes M_n$ of norm 1. Define
\begin{equation*}
 \psi(x)=\sigma(x)+\phi_{(1-p)}(x).
\end{equation*}
Then $\psi$ is a complete order embedding of the required form.  By (\ref{eq:pphippert}) and  (\ref{eq:cornercomp2}) it follows that
\begin{equation*}
 \V \psi-\phi\V_{cb}\leq 57\sqrt{\delta}+\sqrt{115}\delta^{1/4}<68\delta^{1/4}.
\end{equation*}
\end{proof}

\subsection{Non-unital Analogs} 
 Many of our maps of interest will be non-unital CPCs.  For this reason we collect some non-unital analogs (Corollary \ref{cor:nonunitalhomosplit} and Theorem \ref{thm:invertingCCP}) of well-known results about UCP maps that we couldn't find in the literature.

\begin{lemma} Let   $B$ be a unital $C^*$-algebra and  $\phi:M_n\rightarrow B$ a complete order embedding. Let $\tau_n$ denote the normalized trace on $M_n.$  Consider the map $\psi:M_n\rightarrow B$ defined by
\begin{equation}
 \psi(x)=\phi(x)+\tau_n(x)(1-\phi(1)) \label{eq:unitalcoe}
\end{equation}
Then $\psi$ is a (unital) complete order embedding.
\end{lemma}
\begin{proof} $\psi$ is CP as it is the sum of CP maps. To see it is completely isometric, let $x\in (M_k\otimes M_n)^+$, then
 \begin{equation*}
 \V x\V \geq \V\psi^{(k)}(x)\V\geq\V\phi^{(k)}(x)\V=\V x\V.
\end{equation*}
By \cite[Lemma 2.3]{Eckhardt08}, it follows that $\psi$ is a complete isometry.
\end{proof}
We now obtain a non-unital analog of \cite[7.1]{Choi77}.
\begin{corollary} \label{cor:nonunitalhomosplit} Let $B$ be a finite dimensional $C^*$-algebra and $\phi:M_n\rightarrow B$ be a  complete order embedding. Then there is a rank $n$ projection $p\in B$ such that $\phi_p$ is a nonzero *-homomorphism and $\phi=\phi_p+\phi_{(1-p)}.$  
\end{corollary}
\begin{proof} If $n=1,$ this is trivial, so assume $n\geq2.$ Let $\psi$ be as in (\ref{eq:unitalcoe}). By \cite[7.1]{Choi77} there is  a rank $n$ projection $p\in B$ such that $\psi=\psi_p+\psi_{(1-p)}$ with $\psi_p$ a nonzero *-homomorphism.  Let $u\in M_n$ be a unitary with $\tau_n(u)=0.$  Then $\phi(u)=\psi(u)$ and
\begin{equation}
 p\phi(1)p+p(1-\phi(1))p=p\psi(1)p=p\psi(u)p\psi(u)^*p=p\phi(u)p\phi(u)^*p \leq p\phi(1)p. \label{align:cuttingoffstuff}
\end{equation}
Hence, $p(1-\phi(1))p=0$, so $\psi_p=\phi_p.$
\end{proof}
In general we cannot replace $B$ with $B(H)$ for an infinite dimensional Hilbert space in Corollary \ref{cor:nonunitalhomosplit}, because  $\phi(e_{11})$ needn't have 1 as an eigenvalue. But we can get as close as possible (Theorem \ref{thm:homosaredense}).  We start with the following perturbation lemma of Kerr and Li:
\begin{lemma} \textup{(\cite[Lemma 7.1]{Kerr08})} \label{lem:dimdep} For every $\epsilon>0$ and $n\in\mathbb{N}$, there is a $\delta(\epsilon,n)>0$ such that if 
 $\phi:M_n\rightarrow M_N$ is an injective, UCP map with $\V\phi^{-1}\V_{cb}\leq 1+\delta$, then there is a unital complete order embedding $\psi:M_n\rightarrow M_N$ with $\V\phi-\psi\V_{cb}\leq\epsilon.$
\end{lemma}
It is also worth recalling the following result of Roger Smith that we will use repeatedly:
\begin{theorem}[Smith's Lemma] \textup{\cite[Theorem 2.10]{Smith83}} \label{thm:smithslemma} Let $X$ be an operator space and \newline $\phi:X\rightarrow M_n$ a linear map.  Then
\begin{equation*}
 \V \phi^{(n)}\V=\V \phi\V_{cb}.
\end{equation*}
\end{theorem}
\begin{defn}
 Fix $n\in\mathbb{N}$ and a Hilbert space $H.$  Let $\mathcal{H}(n,H)$ denote the set of all complete order embeddings $\phi:M_n\rightarrow B(H)$ such that there is a rank $n$ projection $p\in B(H)$ so $\phi=\phi_p+\phi_{(1-p)}$ and $\phi_p$ is a nonzero *-homomorphism.
\end{defn}

\begin{theorem} \label{thm:homosaredense} $\mathcal{H}(n,H)$ is cb-dense in the set of all complete order embeddings from $M_n$ to $B(H).$
\end{theorem}
\begin{proof} Suppose first that $\phi(1)=1.$ Let $0<\epsilon<(150)^{-1/4}$ and let $\delta=\delta(\epsilon,n)>0$ satisfy Lemma \ref{lem:dimdep}. From Theorem \ref{thm:smithslemma} it follows that
\begin{equation*}
 \V \phi^{-1}:\phi(M_n)\rightarrow M_n\V_{cb}=\V (\phi^{(n)})^{-1}:M_n\otimes \phi(M_n)\rightarrow M_n\otimes M_n\V
\end{equation*}
Since the unit ball of $M_n\otimes \phi(M_n)$ is compact, it follows that there is a finite rank projection $q\in B(H)$ such that 
$\V\phi_q^{-1}\V\leq 1+\delta.$   Use Lemma \ref{lem:dimdep} to obtain a unital complete order embedding $\psi:M_n\rightarrow qB(H)q$
with $\V \phi_q-\psi\V_{cb}\leq \epsilon.$

By \cite[7.1]{Choi77}, there is a rank $n$ projection $p\leq q$ such that $\psi=\psi_p+\psi_{(q-p)}$ and $\psi_p$ is a nonzero *-homomorphism.
Then $\V \phi_p-\psi_p\V_{cb}\leq \epsilon.$  In particular, $\V\phi_p^{-1}\V_{cb}\leq 1+2\epsilon.$  By Lemma \ref{lemma:cutdownalmostiso}, there is a $\widetilde{\psi}\in \mathcal{H}(n,H)$ such that $\V\phi-\widetilde{\psi}\V_{cb}\leq 150\epsilon^{1/4}.$ 

To prove the non-unital case, we first define $\psi$ as in (\ref{eq:unitalcoe}).  Then apply the first part of the proof and use the ``approximate'' version of (\ref{align:cuttingoffstuff}) to obtain the conclusion.
\end{proof}
In the unital case, the following theorem is an immediate consequence of Arveson's extension theorem. The reason for the following theorem is that we don't know if $1$ is in the range of $\phi$, i.e. $\phi(M_n)$ needn't be a sub operator system of $B(H).$  
\begin{theorem} \label{thm:invertingCCP} Let $n\in\mathbb{N}$ and $\phi:M_n\rightarrow B(H)$ be a complete order embedding.  Then there is a UCP map $T:B(H)\rightarrow M_n$
such that $T\phi=id_{M_n}.$ 
\end{theorem}
\begin{proof}  Let $(\psi^k)_{k=1}^\infty$ be a sequence from $\mathcal{H}(n,H)$ that converges to $\phi.$  Let $(p_k)_{k=1}^\infty$ be a sequence of rank $n$ projections such that 
\begin{equation*}
 \psi^k=\psi_{p_k}^k+\psi_{(1-p_k)}^k\quad \textrm{ with }\quad\psi_{p_k}^k\quad\textrm{ a non-zero *-homomorphism.}
\end{equation*}
For each $k\in\mathbb{N}$ define $T_k:B(H)\rightarrow M_n$ by 
\begin{equation}
 T_k(x)=(\psi_{p_k}^k)^{-1}(p_kxp_k).  \label{eq:T_kdefinition}
\end{equation}
Then each $T_k$ is UCP and $T_k\psi^k=id_{M_n},$ hence
\begin{equation*}
\V T_k\phi -id_{M_n}\V_{cb}\leq \V \phi-\psi^k\V_{cb}. 
\end{equation*}
  Finally, let $\omega$ be a non-principal ultrafilter on $\mathbb{N}.$  Then the map
\begin{equation*}
 T(x)=\lim_{k\rightarrow\omega} T_k(x)
\end{equation*}
has the desired properties.
\end{proof}
\subsection{The Case where dim$(H)>n$} \label{sec:varyingrange}
We now turn our attention to the general case,  where the range of $\phi$ is $B(H)$ for an arbitrary Hilbert space.
We start with an example showing that the dependence on $n$ in Lemma \ref{lem:dimdep} cannot be omitted.
\subsubsection{Dependence on Dimension} \label{sec:dimensiondep}
Fix $n>4.$ Let
\begin{equation}
 X_n=\{p\in M_n: p\textrm{ is a rank }n-1\textrm{ projection}\}. \label{eq:X_ndef}
\end{equation}
Let 
\begin{equation}
p_1,...,p_r\in X_n\quad \textrm{ be a }\quad \frac{1}{n}\textrm{-net}\quad \textrm{ for }\quad X_n.\label{eq:finitenetforX_n}
\end{equation} 
Define $\phi:M_r\rightarrow M_r\otimes M_n$ by
\begin{equation*}
 \phi(e_{ij})=e_{ij}\otimes p_ip_j,\quad\textrm{ for }1\leq i,j\leq r.
\end{equation*}
\begin{lemma}
 $\phi$ is a CPC.
\end{lemma}
\begin{proof}
\begin{align*}
 \phi^{(r)}\Big(\sum_{i,j=1}^r e_{ij}\otimes e_{ij}\Big)&=\sum_{i,j=1}^r e_{ij}\otimes e_{ij}\otimes p_ip_j\\
&=\Big(\sum_{i=1}^r e_{i1}\otimes e_{i1}\otimes p_i\Big)\Big(\sum_{i=1}^r e_{i1}\otimes e_{i1}\otimes p_i\Big)^*\geq0.
\end{align*}
Then $\phi$ is completely positive by \cite[Theorem 3.14]{Paulsen02}.  Furthermore, $\phi(1)$ is a projection, hence $\phi$ is completely contractive.
\end{proof}
\begin{lemma}
 $\phi$ is injective with $\V\phi^{-1}\V_{cb}\leq \frac{n}{n-4}.$
\end{lemma}
\begin{proof}
Let $k\in\mathbb{N}$ be arbitrary. Let $p\in M_k\otimes M_r$ be a rank 1 projection.  Then there exist scalars $\alpha_{ij}\in\mathbb{C}$
with $1\leq i\leq k$ and $1\leq j\leq r$ and $\sum_{i,j}|\alpha_{ij}|^2=1$ and an operator $v\in M_k\otimes M_r$ with $vv^*=p$ and
\begin{equation*}
 v=\sum_{i=1}^k e_{i1}\otimes\Big(\sum_{j=1}^r \alpha_{ij}e_{j1}\Big).
\end{equation*}
Let $\tau_n$ denote the normalized trace on $M_n.$ Since each $p_i\in M_n$ is a rank $n-1$ projection, it follows that 
\begin{equation*}
 \tau_n(p_1p_ip_1)=\tau_n(p_i)-\tau_n((1-p_1)p_i)\geq \frac{n-1}{n}-\V p_i\V \tau_n(1-p_i)=\frac{n-2}{n},\quad \textrm{for }i=1,...,r.
\end{equation*}

Therefore
\begin{align*}
 \V \phi^{(k)}(p)\V&\geq \V\phi^{(k)}(v)\phi^{(k)}(v)^*\V\\
&=\V\phi^{(k)}(v)^*\phi^{(k)}(v)\V\\
&=\Big\V \Big(\sum_{i=1}^k e_{1i}\otimes\Big(\sum_{j=1}^r \overline{\alpha}_{ij}e_{1j}\otimes p_1p_j\Big)\Big)\Big(\sum_{i=1}^k e_{i1}\otimes\Big(\sum_{j=1}^r \alpha_{ij}e_{j1}\otimes p_jp_1\Big)\Big)\Big\V\\
&=\Big\V e_{11}\otimes e_{11}\otimes\Big(\sum_{i=1}^k\sum_{j=1}^r|\alpha_{ij}|^2p_1p_jp_1\Big)\Big\V\\
&\geq \tau_n\Big(\sum_{i=1}^k\sum_{j=1}^r|\alpha_{ij}|^2p_1p_jp_1\Big)\\
&\geq \frac{n-2}{n}.
\end{align*}
Now, let $a\in M_k\otimes M_r$ be positive and norm 1.  Then there is a rank 1 projection $p\leq a.$  Since $\phi$ is CP, it follows that
$\V \phi^{(k)}(a)\V\geq \V \phi^{(k)}(p)\V\geq \frac{n-2}{n}.$ Finally, by \cite[Lemma 2.3]{Eckhardt08} it follows that $\phi$ is injective with
\begin{equation*}
 \V \phi^{-1}\V_{cb}\leq \Big( 2\Big(\frac{n-2}{n}\Big)-1\Big)^{-1}=\frac{n}{n-4}.
\end{equation*}
\end{proof}
\begin{theorem}
 Let $\psi:M_r\rightarrow M_r\otimes M_n$ be a complete order embedding.  Then $\V\psi-\phi\V\geq 1-\frac{1}{n}.$
\end{theorem}
\begin{proof}
Since $\psi$ is a complete order embedding, Corollary \ref{cor:nonunitalhomosplit}  provides a rank $r$ projection $p\in M_r\otimes M_n$ such that $\psi_p$ is a *-monomorphism.
In particular, there is a norm 1 vector $\xi\in \ell^2(r)\otimes \ell^2(n)$ such that 
\begin{equation}
 \V \psi(e_{i1})\xi\V=1\quad \textrm{ for every }i=1,...,r. \label{eq:psigood}
\end{equation}
Decompose $\xi=e_1\otimes \xi_1+\eta$ where $\eta\perp e_1\otimes \ell^2(n).$ Let $q\in X_n$ (see Definition (\ref{eq:X_ndef})) be the rank $n-1$ projection onto $(\mathbb{C}p_1(\xi_1))^\perp$ (or if $p_1(\xi_1)=0$, let $q$ be any element of $X_n$).  Then $qp_1(\xi_1)=0.$
By (\ref{eq:finitenetforX_n}), there is an $1\leq i'\leq r$ such that
$\V q-p_{i'}\V< 1/n.$ 

Hence,
\begin{equation*}
 \V \phi(e_{i'1})\xi\V= \V (e_{i'1}\otimes p_{i'}p_1)\xi\V=\V e_{i'}\otimes (p_{i'}p_1\xi_1)\V< \V e_{i'}\otimes (qp_1\xi_1)\V+1/n=\frac{1}{n}. 
\end{equation*}
Combining this with (\ref{eq:psigood}), it follows that 
\begin{equation*}
 \V \psi-\phi\V\geq \V (\psi-\phi)(e_{i'1})\V>1-\frac{1}{n}.
\end{equation*}
\end{proof}
\begin{remark}
 Note that $\phi(1)$ is a projection.  Hence we may additionally assume that $\phi$ is unital above.
\end{remark}

\begin{corollary} \label{cor:counterex} For every $\epsilon>0$ there are $n,N\in \mathbb{N}$ and an injective UCP map $\phi:M_n\rightarrow M_N$ with 
$\V\phi^{-1}\V_{cb}\leq 1+\epsilon$ and $\V\phi-\psi\V\geq 3/4$ for every complete order embedding $\psi.$ 
\end{corollary}
\subsubsection{Independence of Dimension in Special cases}

Let $\phi:M_n\rightarrow B(H)$ be a CPC with $\V \phi^{-1}\V_{cb}\sim 1.$
Corollary \ref{cor:counterex} shows that we cannot perturb $\phi$ to a complete order embedding unless we take dimension into account.
In this section, we show that in two special cases we can perturb $\phi$ independent of dimension.  In particular, we show if the range of $\phi$ is faithful modulo the compact operators (Theorem \ref{thm:antiliminalperturb}) then we can perturb independent of dimension. Next we show (Theorem \ref{thm:amplificationperturb}) that we can always perturb an amplification of $\phi$ regardless of dimension.
Finally, we use Theorem \ref{thm:amplificationperturb} to show that if the range of $\phi$ approximately contains a minimal projection then we can perturb $\phi$ regardless of dimension.  The key to all three results is Lemma \ref{lem:cruxlemma}.
\begin{defn}  Let $n\in\mathbb{N}.$  Let $\textrm{Tr}$ denote the (non-normalized) trace on $M_n.$  Let
\begin{equation}
 \Lambda_n=\{x\in M_n^+:\textrm{Tr}(x)=1\} \label{eq:LAMBDA_n}.
\end{equation} 
For a linear map $\phi:M_n\rightarrow B(H),$ let 
\begin{equation}
 \Lambda_n^\phi=\Big\{\sum_{i,j=1}^n \lambda_{ij} \phi(e_{1i})\phi(e_{j1}):(\lambda_{ij})\in \Lambda_n\Big\}. \label{eq:lambdanphi}
\end{equation}
\end{defn}
\begin{lemma} \label{lem:describesconvexLphiinBH} Let $\phi:M_n\rightarrow B(H)$ be an injective, CPC with $\V\phi^{-1}\V_{cb}\leq 1+\delta.$  Then $\Lambda_n^\phi$ is a compact, convex subset of $B(H)^+$ with
\begin{equation}
x\leq \phi(e_{11})\quad \textrm{ and }\quad \V x\V > \frac{1}{1+3\delta}\quad \textrm{ for every }x\in \Lambda_n^\phi. \label{eq:normestimateonLphi}
\end{equation} 
\end{lemma}
\begin{proof} It is clear that $\Lambda_n$ is compact and convex, which implies that $\Lambda_n^\phi$ shares the same properties. 

Let $\lambda=(\lambda_{ij})_{i,j=1}^n \in \Lambda_n.$ By spectral theory,
$\lambda=\sum_{\ell=1}^k r_\ell p_\ell$ with $r_\ell\in\mathbb{R}^+$ and $p_\ell$ orthogonal, rank 1 projections.   In particular, for each $\ell=1,...,k$ there are scalars $\alpha_1(\ell), ..., \alpha_n(\ell)\in \mathbb{C}$ such that
\begin{equation}
 \lambda_{ij}=\sum_{\ell=1}^k \alpha_i(\ell)\overline{\alpha_j(\ell)}.   \label{eq:lambdaassum}
\end{equation}
Define
\begin{equation*}
v=\left[ \begin{array}{ccccc} \sum_{i=1}^n \alpha_i(1)e_{1i} & \sum_{i=1}^n \alpha_i(2)e_{1i} & \cdot&\cdot& \sum_{i=1}^n \alpha_i(k)e_{1i} \end{array} \right] \in M_{1, k}(M_n).
\end{equation*}
 Then by (\ref{eq:lambdaassum}),
\begin{align}
 e_{11}\otimes\Big(\sum_{i,j=1}^n \lambda_{ij}\phi(e_{1i})\phi(e_{j1})\Big)&=\phi^{(k)}(v)\phi^{(k)}(v)^*  \notag \\
& \leq \phi^{(k)}(vv^*)=e_{11}\otimes (Tr(\lambda))\phi(e_{11})= e_{11}\otimes \phi(e_{11}). \label{eq:descriptionofxinlambda}
\end{align}
This shows that $\Lambda_n^\phi\subset B(H)^+$ and also proves the first inequality from (\ref{eq:normestimateonLphi}). For the second inequality, first notice that
 $\V v\V=1.$
Combining this fact with  (\ref{eq:descriptionofxinlambda}) we have,
\begin{align*}
 \Big\V \sum_{i,j=1}^n \lambda_{ij}\phi(e_{1i})\phi(e_{j1})\Big\V&= \V \phi^{(k)}(v)\phi^{(k)}(v)^*\V \\
&=\V \phi^{(k)}(v)\V^2\geq (1+\delta)^{-2}> (1+3\delta)^{-1}.
\end{align*}
This proves the second inequality from (\ref{eq:normestimateonLphi}).
\end{proof}

\begin{lemma} \label{lem:cruxlemma}Let $\phi:M_n\rightarrow B(H)$ be a CPC and $\delta>0$. Suppose there is a norm 1 vector $\xi\in H$ such that
\begin{equation}
 \min_{x\in\Lambda_n^\phi}\la x\xi,\xi\ra\geq \frac{1}{1+\delta}. \label{eq:minimalvector}
\end{equation}
Then there is a complete order embedding $\psi:M_n\rightarrow B(H)$ such that $\V \psi-\phi\V_{cb}\leq 136\delta^{1/4}.$
 
\end{lemma}
\begin{proof} We will construct a rank $n$ projection $p\in B(H)$ that satisfies  Lemma \ref{lemma:cutdownalmostiso}.

Since $\phi$ is contractive and by (\ref{eq:minimalvector}), for each $\alpha_1, ..., \alpha_n\in \mathbb{C}$ with $\sum_i|\alpha_i|^2=1$ we have
\begin{equation}
 1\geq \Big\V \sum_{i=1}^n \alpha_i \phi(e_{i1})\xi\Big\V=\Big\la \sum_{i,j=1}^n \overline{\alpha}_i\alpha_j\phi(e_{1i})\phi(e_{j1})\xi,\xi\Big\ra^{1/2}\geq \frac{1}{(1+\delta)^{1/2}}. \label{eq:basischange}
\end{equation}
In particular, the vectors $\phi(e_{11})\xi, ..., \phi(e_{n1})\xi$ are linearly independent, so the subspace
\begin{equation*}
 K=\textrm{span}\{ \phi(e_{11})\xi, ..., \phi(e_{n1})\xi\}\subseteq H
\end{equation*}
is $n$ dimensional. Let $p$ be the orthogonal projection from $H$ onto $K.$

Let $q\in M_2\otimes M_n$ be a rank 1 projection.  Then there are scalars $\alpha_1, ..., \alpha_n, \beta_1, ..., \beta_n\in \mathbb{C}$
with $\sum_i |\alpha_i|^2+|\beta_i|^2=1$ such that $q=vv^*$ where
\begin{equation*}
 v=e_{11}\otimes \Big(\sum_{i=1}^n \alpha_i e_{i1}\Big)+e_{21}\otimes\Big(\sum_{i=1}^n \beta_ie_{i1}\Big).
\end{equation*}
Set 
\begin{equation*}
 \eta= e_1\otimes\Big(\sum_{i=1}^n \alpha_i \phi(e_{i1})\xi\Big)+  e_2\otimes\Big(\sum_{i=1}^n \beta_i \phi(e_{i1})\xi\Big)\in \ell^2(2)\otimes K.
\end{equation*}
Then $\V\eta\V\leq1$ by (\ref{eq:basischange}). Furthermore,
\begin{equation}
\la \phi_p^{(2)}(q)\eta,\eta\ra= \la \phi^{(2)}(q)\eta,\eta\ra\geq \la \phi^{(2)}(v)\phi^{(2)}(v)^*\eta,\eta\ra=\V \phi^{(2)}(v)^*\eta\V^2. \label{eq:normpreservedonrank1projections}
\end{equation}
Note that $(\overline{\alpha}_i\alpha_j+\overline{\beta}_i\beta_j)_{i,j}\in \Lambda_n.$ So it follows from (\ref{eq:minimalvector}) and the Cauchy-Schwarz inequality that 

\begin{equation*}
 \V \phi^{(2)}(v)^*\eta\V=\Big\V \sum_{i,j=1}^n (\overline{\alpha}_i\alpha_j+\overline{\beta}_i\beta_j)\phi(e_{1i})\phi(e_{j1})\xi\Big\V\geq  \frac{1}{1+\delta}.
\end{equation*}
Combining this with (\ref{eq:normpreservedonrank1projections}), it follows that
\begin{equation*}
 \V \phi_p^{(2)}(q)\V \geq (1+\delta)^{-2}\quad \textrm{ for every rank 1 projection}\quad q\in M_2\otimes M_n.
\end{equation*}
Let $a\in M_2\otimes M_n$ be positive and norm 1. Then there is a rank 1 projection $q\leq a.$  Since $\phi_p$ is completely positive, it follows that  
\begin{equation*}
\V \phi_p^{(2)}(a)\V\geq \V \phi_p^{(2)}(q)\V \geq \frac{1}{(1+\delta)^2}. 
\end{equation*}
Applying \cite[Lemma 2.3]{Eckhardt08}, it follows that $\phi_p$ is injective with 
\begin{equation*}
 \V \phi_p^{-1}\V \leq \Big( \frac{2}{(1+\delta)^2}-1\Big)^{-1}\leq 1+7\delta.
\end{equation*}
So, by Lemma \ref{lemma:cutdownalmostiso} there is a complete order embedding $\psi:M_n\rightarrow B(H)$ such that
$\V \phi-\psi\V_{cb}\leq 136 \delta^{1/4}.$
\end{proof}

Since von Neumann's minimax theorem appeared in \cite{vonNeumann28}, numerous generalizations have followed.  We will use one such generalization due to Ky Fan in the following lemma. 

\begin{lemma} \label{lem:minimax} Let $A$ be a $C^*$-algebra and $X\subset A^+$ compact and convex such that 
\begin{equation*}
 \V x\V \geq r\quad \textrm{ for all }x\in X.
\end{equation*}
Then there is a state $\omega\in A^*$ such that
\begin{equation*}
  \omega(x) \geq r\quad \textrm{ for all }x\in X.
\end{equation*}
If $A$ is a von Neumann algebra, then $\omega$ can be chosen to be normal.
\end{lemma}
\begin{proof} Let $S(A)$ denote the state space of $A$ equipped with the  $\sigma(A^*, A)$-topology. Then $S(A)$ is compact and convex.
Consider the mapping 
\begin{equation*}
 f:S(A)\times X\rightarrow \mathbb{R}^+\quad \textrm{ defined by }f(\omega, x)=\omega(x).
\end{equation*}
Then $f$ is continuous and affine in each variable.  By \cite[Theorem 3]{Fan52} we have
\begin{equation*}
 \max_{\omega\in S(A)}\min_{x\in X} f(\omega, x)=\min_{x\in X}\max_{\omega\in S(A)} f(\omega, x)\geq r.
\end{equation*}
For the von Neumann case we can replace $S(A)$ with the space of normal states on $A.$
\end{proof}

Let us pause for a moment and map out the rest of the section. Let $\phi:M_n\rightarrow B(H)$ be a CPC with $\V\phi^{-1}\V_{cb}\sim 1.$ 
Then Lemmas \ref{lem:describesconvexLphiinBH} and \ref{lem:minimax} provide a state $\omega$ such that $\omega(x)\sim1$ for $x\in \Lambda_n^\phi.$ Lemma \ref{lem:cruxlemma} then says we can perturb $\phi$ if $\omega$ can be chosen to be a vector state.  For example, for the ``non-perturbable'' map $\phi$ defined in section \ref{sec:dimensiondep} we had $\omega=\la (\cdot)\textrm{ } e_1,e_1\ra\otimes\tau_n$, which is far away from a vector state.  We finish this section with two extreme cases where $\omega$ can be chosen to be a vector state.

\begin{theorem} \label{thm:antiliminalperturb}Let $\phi:M_n\rightarrow B(H)$ be a CPC with $\V \phi^{-1}\V_{cb}<1+\delta.$
Moreover, assume that $C^*(\phi(M_n))\cap K(H)=\{0\}.$  Then there is a complete order embedding $\psi:M_n\rightarrow B(H)$
such that $\V \phi-\psi\V_{cb}\leq 272 \delta^{1/4}.$ 
 
\end{theorem}
\begin{proof}
 Consider $\Lambda_n^\phi$ as in (\ref{eq:lambdanphi}).
By Lemmas \ref{lem:describesconvexLphiinBH} and \ref{lem:minimax} there is a state $\omega\in C^*(\phi(M_n))$ such that 
\begin{equation*}
 \min_{x\in \Lambda_n^\phi} \omega(x) > \frac{1}{1+3\delta}.
\end{equation*}
Since $C^*(\phi(M_n))\cap K(H)=\{0\},$  Glimm's lemma (see \cite[Lemma 11.2.1]{Dixmier77}) states that $\omega$ is a weak* limit of vector states.  Since $\Lambda_n^\phi$ is compact, we obtain a unit vector $\xi\in H$
such that
\begin{equation*}
 \min_{x\in \Lambda_n^\phi}\la x\xi,\xi\ra\geq \frac{1}{1+3\delta}. 
\end{equation*}
The conclusion now follows from Lemma \ref{lem:cruxlemma}.
\end{proof}
Before we consider the second extreme case, we first need to show that we can always perturb amplifications of injective, CPC maps.  We start with a technical lemma.

\begin{lemma} \label{lem:usefulcalculations} Let $B$ be a  $C^*$-algebra and $n\in\mathbb{N}.$  Suppose that 
\newline
$\phi:M_n\oplus B\rightarrow B(H)$ is a CPC with $\V\phi^{-1}\V_{cb}<1+\delta$ and there is a complete order embedding $\psi:M_n\rightarrow B(H)$ such that $\V \phi|_{M_n}-\psi\V_{cb}\leq\delta.$ Moreover assume that there is a rank $n$ projection $q\in B(H)$ such that $\psi=\psi_q+\psi_{(1-q)}$  with $\psi_q$ a nonzero *-homomorphism.   Then, 
\begin{equation}
 \V (\phi_q)|_{0\oplus B}\V_{cb} \leq 8\delta,  \label{eq:Lshapednormcontrol}
\end{equation}
and 
\begin{equation}
\V [\phi(x), q]\V\leq 8\sqrt{\delta}\quad \textrm{ for all }\quad x\in M_n\oplus B \quad \textrm{of norm 1.}\label{eq:cutdownalmostcommute}
\end{equation}
\end{lemma}
\begin{proof} Let $k\in\mathbb{N}$ and  $x\in 0\oplus (M_k\otimes B)$ be positive and norm 1. Let $p=1_n\oplus0\in M_n\oplus B.$  Set $\widetilde{q}:=1_k\otimes q.$ Since $\widetilde{q}\psi^{(k)}(1_k\otimes p)=\widetilde{q}$ and $\V \psi-\phi|_{M_n}\V_{cb}\leq\delta$ we have
\begin{align*}
1&\geq \V \widetilde{q}\phi^{(k)}(1_k\otimes p+x)\V\\
&\geq \V \widetilde{q}+\widetilde{q}\phi^{(k)}(x)\V-\delta\\
&= \V \widetilde{q}+\widetilde{q}\phi^{(k)}(x)\widetilde{q}+\widetilde{q}\phi^{(k)}(x)(1-\widetilde{q})\V-\delta\\
&=\V \widetilde{q}+2\widetilde{q}\phi^{(k)}(x)\widetilde{q}+\widetilde{q}\phi^{(k)}(x)\widetilde{q}\phi^{(k)}(x)\widetilde{q}+\widetilde{q}\phi^{(k)}(x)(1-\widetilde{q})\phi^{(k)}(x)\widetilde{q}\V^{1/2}-\delta\\
&=\Big(1+\V2\widetilde{q}\phi^{(k)}(x)\widetilde{q}+\widetilde{q}\phi^{(k)}(x)\widetilde{q}\phi^{(k)}(x)\widetilde{q}+\widetilde{q}\phi^{(k)}(x)(1-\widetilde{q})\phi^{(k)}(x)\widetilde{q}\V\Big)^{1/2}-\delta.
\end{align*}
After rearranging it follows that
\begin{equation}
\max\{ \V 2\widetilde{q}\phi^{(k)}(x)\widetilde{q}\V, \V \widetilde{q}\phi^{(k)}(x)(1-\widetilde{q})\phi^{(k)}(x)\widetilde{q}\V\}\leq 2\delta+\delta^2. \label{eq:2for1inequality}
\end{equation}
Since every $y\in 0\oplus(M_k\otimes B)$ of norm 1 can be written as the sum of four positive elements each with norm bound by 1, we have
\begin{equation*}
 \V (\phi_q)|_{0\oplus B}\V_{cb}\leq 4(\delta+\delta^2/2)\leq 8\delta.
\end{equation*}
Let $(z,x)\in M_n\oplus B$ be positive and norm 1.  By hypothesis, we have $\V[\phi(z,0),q]\V\leq \delta.$ By (\ref{eq:2for1inequality}), we have $\V[\phi(0,x),q]\V\leq (2\delta+\delta^2)^{1/2}.$  Again by decomposing an arbitrary $y\in M_n\oplus B$ into positive pieces, we obtain 
(\ref{eq:cutdownalmostcommute}).
\end{proof}

\begin{defn} Let $\phi:X\rightarrow Y$ be a linear map and  $k\in\mathbb{N}.$
Define $1_k\otimes \phi:X\rightarrow M_k(Y)$ by
\begin{equation*}
(1_k\otimes \phi)(x)=1_k\otimes\phi(x).
\end{equation*} 
\end{defn}

\begin{theorem} \label{thm:amplificationperturb} Let $n\in\mathbb{N}$  and let $\phi:M_n\rightarrow B(H)$ be a CPC with $\V \phi^{-1}\V_{cb}<1+\delta.$ Then there is a natural number $k\in\mathbb{N}$ and a complete order embedding $\psi:M_n\rightarrow M_k\otimes B(H)$ such that
\begin{equation*}
 \V 1_k\otimes\phi-\psi\V_{cb}\leq 272 \delta^{1/4}.
\end{equation*}
 \end{theorem}
\begin{proof} 
By Lemmas \ref{lem:describesconvexLphiinBH} and \ref{lem:minimax} we obtain a state $\omega\in B(H)_*$ such that
\begin{equation*}
 \min_{x\in \Lambda_n^\phi}\omega(x) > \frac{1}{1+3\delta}. 
\end{equation*}
Since $\omega\in B(H)_*$  there are  vectors $\xi_1,\xi_2,...\in H$ such that
\begin{equation*}
 \omega(x)=\sum_{i=1}^\infty \la x \xi_i,\xi_i\ra\quad \textrm{ for all }x\in B(H).
\end{equation*}
Since $\Lambda_n^\phi$ is compact, there is a $k\in\mathbb{N}$ such that
\begin{equation}
 \min_{x\in \Lambda_n^\phi}\sum_{i=1}^k\la x\xi_i,\xi_i\ra \geq \frac{1}{1+3\delta}.  \label{eq:minimalstate}
\end{equation}

Let $\xi=\sum_{i=1}^k e_i\otimes \xi_i\in H^k.$  Then for any $x\in B(H)$ we have
\begin{equation*}
 \la (1_k\otimes x) \xi, \xi\ra=\sum_{i=1}^k\la x\xi_i,\xi_i\ra.
\end{equation*}
Combining this with (\ref{eq:minimalstate}) and Lemma \ref{lem:cruxlemma} the conclusion follows.

\end{proof}

If we assume that $\phi$ is unital in Proposition \ref{prop:invertingcpc}, the conclusion is an easy consequence of \cite[Proposition 1.19]{Wassermann94} and Wittstock's extension theorem, with a much better estimate.
\begin{proposition} \label{prop:invertingcpc} Let $n\in\mathbb{N}$ and let $\phi:M_n\rightarrow B(H)$ be an injective CPC with $\V \phi^{-1}\V_{cb}<1+\delta.$
Then there is a UCP map, $T:B(H)\rightarrow M_n$ such that 
\begin{equation*}
 \V id_{M_n}-T\phi\V_{cb}\leq 272 \delta^{1/4}.
\end{equation*} 
\end{proposition}
\begin{proof}  By Theorem \ref{thm:amplificationperturb} there is a $k\in\mathbb{N}$ and a complete order embedding $\psi:M_n\rightarrow M_k\otimes B(H)$ such that
\begin{equation*}
 \V 1_k\otimes\phi-\psi\V_{cb}\leq 272 \delta^{1/4}.
\end{equation*}
By Theorem \ref{thm:invertingCCP} there is a UCP map  $R:M_k\otimes B(H)\rightarrow M_n$ such that $R\psi=id_{M_n}.$ Define $T:B(H)\rightarrow M_n$ by $T(x)=R\circ(1_k\otimes x).$  Then,
\begin{equation*}
 \V id_B-T\phi\V_{cb}=\V R\psi-R\circ (1_k\otimes\phi)\V_{cb}\leq \V \psi- 1_k\otimes\phi\V_{cb}\leq 272 \delta^{1/4}.
\end{equation*}

\end{proof}
We are now ready to prove our perturbation theorem for the second extreme case.
\begin{theorem} \label{thm:minimalperturb} Let $\phi:M_n\rightarrow B(H)$ be a CPC with $\V \phi^{-1}\V_{cb}<1+\delta.$
Suppose there is a rank 1 projection $p\in B(H)$ and  $x\in M_n$ such that $\V \phi(x)-p\V <\delta.$  
Then there are
\begin{enumerate}
 \item a rank 1 projection $r\in M_n$ such that $\V \phi(r)-p\V\leq 315\delta^{1/8}$
\item a complete order embedding $\psi:M_n\rightarrow B(H)$ such that $\V \psi-\phi\V_{cb}\leq 1360\delta^{\frac{1}{32}}.$
\end{enumerate}
\end{theorem}
\begin{proof} We first prove (1). 

 Let $T:B(H)\rightarrow M_n$ be as in Proposition \ref{prop:invertingcpc}.  Then $T(p)\geq0$ and 
\begin{align*}
\V \phi(T(p))-p\V &\leq \V \phi(T(p))-\phi(x)\V+\delta\\
&\leq\V T(p)-x\V+\delta\\
&\leq \V T(p)-T(\phi(x))\V+\V T(\phi(x))-x\V+\delta\\
&\leq \V p-\phi(x)\V+272 \delta^{1/4}+\delta\leq 273\delta^{1/4}.
\end{align*}
Then, $\V T(p)\V\geq \V \phi(T(p))\V>1-273\delta^{1/4}.$  
\begin{equation}
\textrm{Let }q=\V T(p)\V^{-1}T(p),\quad \textrm{ then }\quad\V \phi(q)-p\V\leq 819\delta^{1/4}. \label{eq:qclosetop} 
\end{equation}
Let $\delta'=819\delta^{1/4}.$  Let $r\leq q$ be a rank 1 projection. Let $\xi\in H$ be norm 1 such that $p\xi=\xi.$
\newline
Let $\eta\in H$ of norm 1 with $\la \phi(r)\eta,\eta\ra\geq (1+\delta)^{-1}.$ We will now show:
\begin{equation}
 (\exists \theta\in \mathbb{R})(\V \eta-e^{i\theta}\xi\V \leq 2\sqrt{\delta'})  \label{eq:claimformin}
\end{equation}
By (\ref{eq:qclosetop}), it follows that
\begin{equation*}
 \delta'\geq \la (\phi(r)+\phi(q-r)-p)\eta,\eta\ra\geq \la \phi(r)\eta,\eta\ra-\la p\eta,\eta\ra\geq (1+\delta)^{-1}-\la p\eta,\eta\ra.
\end{equation*}
Hence,
\begin{equation*}
 |\la \xi,\eta\ra|=\la p\eta,\eta\ra^{1/2}\geq \Big(\frac{1}{1+\delta}-\delta'\Big)^{1/2}\geq 1-2\delta'.
\end{equation*}
Choose $\theta$ so $\la e^{i\theta}\xi,\eta\ra=|\la \xi,\eta\ra|.$  Then,
\begin{equation*}
 \V e^{i\theta}\xi-\eta\V^2 =2-2|\la \xi,\eta\ra|\leq 4\delta'.
\end{equation*}
This proves (\ref{eq:claimformin}).

Now, suppose that $\xi^\perp$ is norm 1 and perpendicular to $\xi.$  From (\ref{eq:qclosetop}) we deduce that
\begin{equation}
\la \phi(r)\xi^\perp,\xi^\perp\ra \leq \delta'. \label{eq:xiperpest}
\end{equation}
Without loss of generality, suppose that $\phi(r)$ has a $\V\phi(r)\V$ eigenvalue with eigenvector $\zeta.$
By, (\ref{eq:claimformin}) it follows that
\begin{equation}
 \V\phi(r)\xi-\xi\V \leq \V\phi(r)\zeta-\zeta\V+4\sqrt{\delta'}\leq 5\sqrt{\delta'}. \label{eq:xiest}
\end{equation}
Let $t$ be the projection onto $\zeta.$  Then, $\V t-p\V\leq 2\sqrt{\delta'}$ by (\ref{eq:claimformin}).  Thereofore,
\begin{align}
 \V (1-p)\phi(r) p\V &\leq \V (1-t)\phi(r)t\V+2\V t-p\V \notag \\
&\leq \sqrt{1-\V\phi(r)\V}+4\sqrt{\delta'}\leq 5\delta'. \label{align:xicornerest}
\end{align}
Combining (\ref{eq:xiperpest}), (\ref{eq:xiest}) and (\ref{align:xicornerest}) it follows that
\begin{equation*}
 \V \phi(r)-p\V\leq 11\sqrt{\delta'}\leq 315\delta^{1/8}.
\end{equation*}

We now prove (2).

Note that for any unitary $u\in M_n$ the existence of a perturbation for $\phi$ is equivalent to the existence of a perturbation for the map $x\mapsto \phi(uxu^*).$  So without loss of generality, assume that $r=e_{11}.$
By Lemma \ref{lem:describesconvexLphiinBH}, we have 
\begin{equation*}
 \V x\V\geq (1+\delta)^{-1}\quad \textrm{ and }\quad x\leq \phi(e_{11})\quad\textrm{ for every }x\in \Lambda_n^\phi.
\end{equation*}
Let $x\in\Lambda_n^\phi$ and $\eta\in H$ norm 1 such that $\la x\eta,\eta\ra\geq (1+\delta)^{-1}.$  Then
\begin{equation*}
 \la \phi(e_{11})\eta,\eta\ra\geq \la x\eta,\eta\ra\geq (1+\delta)^{-1}.
\end{equation*}
So, by (\ref{eq:claimformin}) we have
\begin{equation*}
 \la \phi(x)\xi,\xi\ra \geq (1+\delta)^{-1}-3(2\sqrt{\delta'})\geq \frac{1}{1+12\sqrt{\delta'}} \quad \textrm{for every }x\in \Lambda_n^\phi .
\end{equation*}
Finally, by Lemma \ref{lem:cruxlemma} there is a complete order embedding $\psi:M_n\rightarrow B(H)$ such that
\begin{equation*}
 \V \psi-\phi\V_{cb}\leq 136(12\sqrt{\delta'})^{1/2}\leq 1360 \delta^{\frac{1}{32}}.
\end{equation*} 
\end{proof}
\begin{corollary} \label{cor:liftingrank1projectionsgeneral} Let $B$ be a finite dimensional $C^*$-algebra and $\phi:B\rightarrow B(H)$  a CPC with $\V \phi^{-1}\V_{cb}<1+\delta.$
Suppose there is a rank 1 projection $p\in B(H)$ and a $x\in B$ such that $\V \phi(x)-p\V <\delta.$  
Then there is a rank 1 projection $r\in B$ such that $\V \phi(r)-p\V\leq 315\delta^{1/8}.$
\end{corollary}
\section{Strong NF algebras}
  Before we can apply our perturbation results, we first need some technical theorems.
\subsection{Technical Theorems}
\begin{defn}
 Let $I$ be an index set and $(n_i)_{i\in I}$ a family of positive integers.  Define the $C^*$-algebra
\begin{equation*}
 \mathcal{B}(I,(n_i)):=\prod_{i\in I} M_{n_i}.
\end{equation*}
For a vector space $X$ and a linear map $\phi:X\rightarrow \mathcal{B}(I,(n_i))$, define for each $i\in I$ the linear map 
$\phi_i:X\rightarrow M_{n_i}$ as the $i$th coordinate map.
\end{defn}
Our first goal of this section is to prove the following:
\begin{theorem} \label{thm:spllitingupwhenrangeisarbitrary}
Let $\epsilon>0$ and $n\in\mathbb{N}.$  There is a $\delta=\delta(n,\epsilon)>0$ such that for every $C^*$-algebra $\mathcal{B}=\mathcal{B}(I,(n_i))$ and every injective linear map $\phi:M_n\rightarrow \mathcal{B}$ with $\V\phi^{-1}\V_{cb}<1+\delta$, there is an index $i\in I$ such that $\phi_i$ is injective with $\V\phi_i^{-1}\V_{cb}<1+\epsilon.$
\end{theorem}

In order to prove Theorem \ref{thm:spllitingupwhenrangeisarbitrary} it is necessary to introduce ternary rings of operators (TROs).  If $\phi$ is assumed unital above then we can compose an (essentially identical)  proof using the theory of completely positive maps.  Unfortunately it seems that deducing the non-unital case from the unital case is a non-trivial matter, and we will definitely require the non-unital case for our applications.  For this reason, we now recall some of the basic facts about TROs required for our proof.

Let $A$ be a $C^*$-algebra and $V\subseteq A$ a closed subspace.  If $V$ is closed under the triple product
\begin{equation*}
 (x,y,z)\mapsto xy^*z \quad \textrm{where}\quad x,y,z\in V,
\end{equation*}
 then $V$ is called a \emph{ternary ring of operators} (henceforth TRO).   We refer the reader to \cite[Section 2]{Kaur02} and the references therein for a detailed look at the development of the theory of TROs.  

Consider the following subspaces of $A:$
\begin{equation}
 V^*:=\{ v^*:v\in V\}\quad C(V)=\overline{span}\{vw^*:v,w\in V\}\quad  D(V)=\overline{span}\{v^*w:v,w\in V\}
\end{equation}
It is straightforward to verify that $C(V)$ and $D(V)$ are $C^*$-subalgebras of $A.$  Furthermore, one checks that
\begin{equation*}
 A(V)=\left[ \begin{array}{ll} C(V) & V \\ V^* & D(V)\\ \end{array} \right] \subset M_2\otimes A. 
\end{equation*}
is a $C^*$-algebra.

A linear map $\phi:V\rightarrow W$ between TROs is called a \emph{TRO-homomorphism} if 
\begin{equation*}
 \phi(xy^*z)=\phi(x)\phi(y)^*\phi(z)\quad \textrm{ for all }\quad x,y,z\in V.
\end{equation*}

\begin{theorem}\textup{\cite{Hamana99}} \label{thm:Hamanalinkinghomo}
Let $\phi:V\rightarrow W$ be a TRO homomorphism.  Then there are *-homomorphisms 
\begin{equation*}
 \phi_C:C(V)\rightarrow C(W)\quad \textrm{ and }\quad\phi_D:D(V)\rightarrow D(W)
\end{equation*}
defined by
\begin{equation*}
 \phi_C\Big( \sum v_iw_i^* \Big)=\sum \phi(v_i)\phi(w_i)^* \quad\textrm{ and }\quad \phi_D\Big( \sum v_i^*w_i \Big)=\sum \phi(v_i)^*\phi(w_i).
\end{equation*}
Moreover, the following map
\begin{equation*}
 \pi_\phi:=\left[ \begin{array}{ll} \phi_C & \phi \\ \phi^* & \phi_D\end{array} \right]:A(V)\rightarrow A(W)
\end{equation*}
is a well-defined *-homomorphism.
\end{theorem}
The following theorem of Masamichi Hamana is the key element to our proof:
\begin{theorem} \textup{\cite[Theorem 3.2(ii)]{Hamana99}} \label{thm:Hamanatheorem} Let $V$ and $W$ be TROs and $\phi:V\rightarrow W$ a complete isometry.  Then there is a TRO-homomorphism $T:TRO(\phi(V))\rightarrow V$ extending $\phi^{-1}$, where $TRO(\phi(V))$ is the TRO generated by $\phi(V).$ 
\end{theorem}
We now recall some facts about $n$-homogeneous $C^*$-algebras.
Let $X$ be a compact Hausdorff space and $n\in \mathbb{N}.$  Under the usual identification 
\begin{equation*}
 C(X)\otimes M_n\cong \{ f:X\rightarrow M_n: f \textrm{ is continuous}\},
\end{equation*}
it is easy to see that there is a one-to-one correspondence between elements $x\in X$ and irreducible representations, $\pi_x$ of $C(X)\otimes M_n$ given by
\begin{equation}
 \pi_x(f)=f(x). \label{eq:irreduciblerepsnhomo}
\end{equation}
Now let $I$ be any index set, then
\begin{equation*}
 \prod_{i\in I} M_n \cong \ell^\infty(I)\otimes M_n \cong C(\beta I)\otimes M_n,
\end{equation*}
where $\beta I$ denotes the Stone-Cech compactification of $I.$  Recall that $\beta I$ is identified with the set of all ultrafilters on $I$ (here $I\subset \beta I$ corresponds to the principal ultrafilters on $I$ in the obvious way). Under this identification and by (\ref{eq:irreduciblerepsnhomo}), it follows that there is a one-to-one correspondence between ultrafilters $\omega$ on $I$ and irreducible representations $\pi_\omega$ of   $\prod_{i\in I} M_n$ given by
\begin{equation}
 \pi_\omega((x_i)_{i\in I})=\lim _{i\rightarrow \omega} x_i. \label{eq:irreduciblerepsofsequencespace}
\end{equation}
It is a well-known result of Krein that for $C^*$-algebras $A\subseteq B$, every pure state $\phi$ on $A$ has a pure state extension to $B.$ The following lemma is a consequence of this fact.
\begin{lemma} 
 Let $A\subseteq B$ be $C^*$-algebras.  Let $(\pi,H)$ be an irreducible representation of $A.$  Then there is an irreducible representation $(\rho, K)$ of $B,$ such that $H\subseteq K$ and letting $p:K\rightarrow H$ be the orthogonal projection, we have
\begin{equation*}
 p\rho(x)p=\pi(x)\quad \textrm{ for every }\quad x\in A.
\end{equation*}
In particular, if $A\subseteq B\cong C(X)\otimes M_n$ for some Hausdorff space $X$, and if $\pi$ is an $n$-dimensional irreducible representation of $A$, then the extension $\rho$ is necessarily $n$-dimensional by \textup{(\ref{eq:irreduciblerepsnhomo})}, hence
\begin{equation}
 \rho(x)=\pi(x)\quad \textrm{ for every }\quad x\in A. \label{eq:nhomogeneousextension}
\end{equation}
\end{lemma}

We now combine all of this background material in the following lemma.
\begin{lemma} \label{lem:coordinatealmostisom} Let $n\in\mathbb{N},$  $I$ an index set, and  $\phi:M_n\rightarrow \ell^\infty(I)\otimes M_n$ be a complete isometry.    For every $\epsilon>0$ there exists an index $i\in I$ such that
$\phi_i$ is an injective complete contraction with $\V\phi_i^{-1}\V_{cb}<1+\epsilon.$
\end{lemma}
\begin{proof}  Let $W\subset \ell^\infty(I)\otimes M_n$ be the TRO generated by $\phi(M_n).$ By Theorem \ref{thm:Hamanatheorem},  there is a TRO homomorphism
\begin{equation}
T:W\rightarrow M_n\quad  \textrm{such that }\quad T=\phi^{-1} \textrm{ on }\phi(M_n). \label{eq:TROextension}
\end{equation}
 Let $\pi_T:A(W)\rightarrow A(M_n)=M_{2n}$ be the *-homomorphism
associated with $T,$ as in Theorem \ref{thm:Hamanalinkinghomo}.     Since $T$ is surjective, so is $\pi_T,$ and hence it is irreducible. We extend $\pi_T$ by (\ref{eq:nhomogeneousextension}) (still call it $\pi_T$)
to an irreducible representation of  $\ell^\infty(I)\otimes M_{2n}.$ By (\ref{eq:irreduciblerepsofsequencespace}), there exists an ultrafilter $\omega$ on $I$ such that $\pi_T=\pi_\omega.$

Now, let $x\in M_n\otimes M_n$ be arbitrary.  Then,
\begin{align}
 \V x \V_{M_n\otimes M_n} &= \Big\V \left[ \begin{array}{ll} 0 & x \\ 0 & 0\\ \end{array} \right]\Big\V_{A(M_n\otimes M_n)}
=\Big\V \left[ \begin{array}{cc} 0 & (\phi^{(n)})^{-1}\circ \phi^{(n)}(x) \\ 0 & 0\\ \end{array} \right]\Big\V_{A(M_n\otimes M_n)} \notag \\
&=\Big\V \left[ \begin{array}{cc} 0 & T^{(n)}\circ \phi^{(n)}(x) \\ 0 & 0\\ \end{array} \right]\Big\V_{A(M_n\otimes M_n)}
=\Big\V \pi_\omega^{(n)}\Big(\left[ \begin{array}{cc} 0 & \phi^{(n)}(x) \\ 0 & 0\\ \end{array} \right]\Big)\Big\V_{A(M_n\otimes M_n)} \notag \\
&= \lim_{i\rightarrow\omega}\V \phi_i^{(n)}(x)\V_{M_n\otimes M_n}. \label{align:normcomputationsultra}
\end{align}
Let $\delta>0$ and let $F$ be a finite $\delta$-net for the unit sphere of $M_n\otimes M_n.$  It follows from (\ref{align:normcomputationsultra}) that
\begin{equation*}
 \bigcap_{x\in F} \{ i\in \omega: \V \phi^{(n)}_i(x)\V\geq 1-\delta\} \in\omega.
\end{equation*}
In particular this set is non-empty.  Let $i'\in I$ be in the above set.  Then for every $x$ in the unit sphere of $M_n\otimes M_n$, it follows that $\V\phi_{i'}^{(n)}(x)\V\geq 1-2\delta.$   
By Theorem \ref{thm:smithslemma}, it follows that $\V \phi_{i'}^{-1}\V_{cb}\leq (1-2\delta)^{-1}.$  As $\delta$ was arbitrary, this proves the lemma.
\end{proof}
\begin{corollary} \label{cor:normpreservedoncoordinates}
 Let $\epsilon>0$ and $n\in\mathbb{N}.$  There exists a $\delta=\delta(\epsilon,n)>0$ such that for any index set $I$, if 
 $\phi:M_n\rightarrow \ell^\infty(I)\otimes M_n$ is an injective complete contraction such that
$\V\phi^{-1}\V_{cb}<1+\delta,$ then there is an index $ i\in I$ such that the $i$-th coordinate map
$\phi_i:M_n\rightarrow M_n$ is an injective complete contraction with $\V\phi_i\V_{cb}<1+\epsilon.$
\end{corollary}
\begin{proof}
  Suppose not.  Then, there is an $\epsilon>0,$ $n\in\mathbb{N},$ a sequence of index sets $I_k$,
and injective complete contractions
\begin{equation*}
 \phi^k:M_n\rightarrow \ell^\infty(I_k)\otimes M_n
\end{equation*}
such that 
\begin{equation}
\V(\phi^k)^{-1}\V_{cb}<1+ 1/k,\textrm{ but }\V (\phi^k_i)^{-1} \V_{cb}>1+\epsilon \textrm{ for every }k\in\mathbb{N}, \textrm{ and }i\in I_k.  \label{eq:coordinatemapsarefar}
\end{equation}
It follows that the map
\begin{equation*}
 \phi=\bigoplus_{k\in\mathbb{N}}\phi^k:M_n\rightarrow \ell^\infty(\sqcup_{k=1}^\infty I_k)\otimes M_n
\end{equation*}
is a complete isometry.  But by Lemma \ref{lem:coordinatealmostisom}, one of the coordinate maps of $\phi$  must be injective with $\V\phi_i^{-1}\V_{cb}<1+\epsilon.$  This contradicts  (\ref{eq:coordinatemapsarefar}) 

\end{proof}

Now let $m\leq n$ be positive integers.  Set
\begin{equation}
 I(m,n)=\{ p\in M_n: p \textrm{ is a rank }m\textrm{ projection}\}.
\end{equation}
Let $x \in M_m\otimes M_n$ be norm 1.

Let $\xi=(\xi_1,...,\xi_m),\eta=(\eta_1,...,\eta_m)\in \ell^2(m)\otimes \ell^2(n)$ be norm 1 such that $x\xi=\eta.$ Let $p_\xi,p_\eta\in M_n$ denote the projections onto $\textrm{span}\{\xi_1,...,\xi_m\}$ and $\textrm{span}\{\eta_1,...,\eta_m\}$ respectively. It is then clear that $\V x\V =\V (1_m\otimes p_\eta)x(1_m\otimes p_\xi) \V.$  From which we deduce that the map
$P_{m,n}:M_n\rightarrow \ell^\infty(I(m,n)^2)\otimes M_m$  defined by
\begin{equation}
 P_{m,n}(x)=\bigoplus_{(p,q)\in I(m,n)^2} pxq \label{eq:definitionofP_nm}
\end{equation}
is a complete contraction, such that $P_{m,n}^{(m)}$ is an isometry.

We are now in a position to prove Theorem \ref{thm:spllitingupwhenrangeisarbitrary}.

\begin{proof}[Proof of Theorem \ref{thm:spllitingupwhenrangeisarbitrary}] Let $\epsilon>0$ and $n\in\mathbb{N}.$  Set $\delta=\delta(\epsilon,n)$ from Corollary \ref{cor:normpreservedoncoordinates}.  Let $\phi:M_n\rightarrow \mathcal{B}(I,(n_i))$ be an injective complete contraction with $\V\phi^{-1}\V_{cb}<1+\delta.$  Without loss of generality, suppose that $n_i\geq n$ for each $i\in I$ (if not simply embed $M_{n_i}\hookrightarrow M_{n_i}\oplus M_{n-n_i}\subset M_n$).

Now consider the map $P:\mathcal{B}(I,(n_i))\rightarrow \bigoplus_{i\in I}\bigoplus_{(p,q)\in I(n,n_i)^2} M_n$ defined by (recall (\ref{eq:definitionofP_nm}))
\begin{equation*}
 P((x_i)_{i\in I})=\bigoplus_{i\in I}P_{n,n_i}(x_i)
\end{equation*}
Then $P$ is a complete contraction and $P^{(n)}$ is an isometry.  It follows from Theorem \ref{thm:smithslemma} that $P\circ \phi$  is an injective complete contraction with 
$\V (P\circ\phi)^{-1}  \V_{cb}<1+\delta.$  By  Corollary \ref{cor:normpreservedoncoordinates}, there is an index $i\in I$ and rank $n$ projections $p,q\in M_{n_i}$ such that $x\mapsto p\phi_i(x)q$ is an injective complete contraction with inverse cb norm bound by $1+\epsilon.$  Hence $\phi_i$ is injective with $\V\phi_i^{-1}\V_{cb}<1+\epsilon.$ 

\end{proof}
The following lemma is a variation of \cite[Lemma 2.2]{Eckhardt08}
\begin{lemma} \label{cor:usefulUCPOL} Suppose $A$ is a unital $C^*$-algebra with $\mathcal{OL}_\infty(A)=1.$  Then for every finite subset $F\subset A$ and every $\delta>0$ there exists a finite dimensional $C^*$-algebra $B$ and maps
$\phi:B\rightarrow A$, $\psi:A\rightarrow B$ such that
\begin{enumerate}
 \item $\phi$ and $\psi$ are UCP.
\item $dist(x, \psi(B))<\delta$ for all $x\in F.$
\item $\V \psi\phi-id_B\V<\delta.$
\item $\V\psi(xy)-\psi(x)\psi(y)\V\leq \delta\V x\V\V y\V$ for all $x,y\in F.$
\end{enumerate}\end{lemma}
\begin{proof}
 Obtain $B$, $\psi$ and the unital self-adjoint map $\phi$ as in \cite[Lemma 2.2]{Eckhardt08}, moreover assume $\psi$ satisfies (4) by the proof of \cite[Theorem 3.2]{Junge03}.  Since $A$ is nuclear,  there is a matrix algebra $M_n$ and  UCP maps $\alpha:A\rightarrow M_n$ and $\beta:M_n\rightarrow A$
such that $\V \beta\alpha|_{\phi(B)}-id_{\phi(B)}\V_{cb}<\delta.$ By \cite[Proposition 1.19]{Wassermann94} there is a UCP map $T:B\rightarrow M_n$ such that $\V T-\alpha\phi\V_{cb}<\delta.$  Then $\widetilde{\phi}=\beta T$ and $\psi$ are our desired maps.  
\end{proof}

\subsection{Applications}

\begin{defn} \label{defn:GCRideal} Let \textup{Prim}(A) denote the primitive ideal space of the $C^*$-algebra A and let $J\in \textrm{Prim}(A).$  We say $J$ is a \emph{GCR ideal} if there is an irreducible representation $(\pi, H)$ of $A$ with 
$\textrm{ker}(\pi)=J$ and $\pi(A)\cap K(H)\neq\{0\}$ (hence $K(H)\subset \pi(A)$). 
\end{defn}

\begin{theorem} \label{thm:maintheoremin3}
 Let $A$ be a $C^*$-algebra with $\mathcal{OL}_\infty(A)=1.$ Suppose $J_1, ..., J_n$ are primitive ideals of $A$ such that $J_1\cap\cdots\cap J_n=\{0\}.$  Then $A$ is a strong NF algebra. In particular, all primitive $C^*$-algebras with $\mathcal{OL}_\infty(A)=1$ are strong NF.
\end{theorem}
\begin{proof} 
Without loss of generality, suppose that $\{J_1,...,J_n\}$ is a minimal element (with respect to set inclusion) of the set of all finite subsets $F\subset \textrm{Prim}(A)$ with $\textrm{ker}(F)=\{0\}.$  In particular, $J_i\neq J_k$ if $i\neq k$ and
\begin{equation}
\bigcap_{k\neq i}J_k\neq\{0\}\quad \textrm{ for every } i=1,...,n. \label{eq:nontrivialint} 
\end{equation}
We also order the $J_i$'s as
$J_1,...,J_r, J_{r+1}, ..., J_n$ so that $J_1, ..., J_r$ are all GCR ideals (Definition \ref{defn:GCRideal}) and $J_{r+1},..., J_n$ are not GCR ideals. 

Let $(\pi_1, H_1),..., (\pi_n, H_n)$ be irreducible representations of $A$ with $\textrm{ker}(\pi_i)=J_i$ and let 
\begin{equation*}
(\pi, H)=\oplus_i (\pi_i, H_i).
\end{equation*}
Since all of the $J_i$ are different, it follows from \cite[Theorem 3.8.11]{Pedersen79} that
\begin{equation}
 \pi(A)''=\prod _{i=1}^n \pi_i(A)''=\prod _{i=1}^n B(H_i). \label{eq:mutexclusiveideals}
\end{equation}
We will show that $A$ is inner quasidiagonal, which combined with \cite[Theorem 4.5]{Blackadar01} shows that $A$ is a strong NF algebra. By (\ref{eq:mutexclusiveideals}) we must show that for each finite set $F\subset A$ and $\epsilon>0$, there is a finite rank projection  $p\in\prod _{i=1}^n B(H_i)$
such that
\begin{equation}
 \max_{x\in F}\V[x,p]\V\leq\epsilon\quad\textrm{ and }\quad \min_{x\in F}\V p\pi(x)p\V\geq 1-\epsilon. \label{eq:representationversionQD}
\end{equation}

In order to keep the notation within reason, we first prove the following lemma:
\begin{lemma} \label{lem:GCRQDrep} $\pi_i$ is a quasidiagonal representation for each $i=1,...,r$ (recall these are the GCR representations). 
\end{lemma}
\emph{Proof of Lemma \ref{lem:GCRQDrep}}
\newline
Of course, it is enough to prove this for $i=1.$  Set $I=J_2\cap\cdots\cap J_n.$ By (\ref{eq:nontrivialint}), $I\neq\{0\}.$ Since $J_1=\textrm{ker}(\pi_1)$ and $J_1\cap I=\{0\}$, it follows that 
\begin{equation}
 \pi_1|_{I} \quad \textrm{ is a *-monomorphism.}\label{eq:pirestrictedtoI}
\end{equation}

Moreover, since $K(H_1)\subseteq \pi_1(A)$ and $\pi_1(I)$ is a non-zero ideal of $\pi_1(A)$, we have $\pi_1(I)\cap K(H_1)\neq0$ from which it follows that 
\begin{equation}
K(H_1)\subseteq \pi_1(I). \label{eq:Ihascompacts}
\end{equation}
Let $\epsilon>0$ and $\pi_1(x_1),...,\pi_1(x_s)\in \pi_1(A)$ be norm 1 such that $x_1,...,x_s\in A$ are also norm 1.
Let $r\in K(H_1)$ be any finite rank projection and let $p\in K(H_1)$ be a finite rank projection such that
\begin{equation}
 \min_{i=1,...,s}\V p\pi_1(x_i)p\V\geq 1-\epsilon/2\quad \textrm{ and }\quad r\leq p. \label{eq:QDconditions}
\end{equation}
Set $N=\textrm{rank}(p)$ and let $(e_{ij})_{i,j=1}^N$ be matrix units for $pB(H_1)p.$ 
By (\ref{eq:pirestrictedtoI}) and (\ref{eq:Ihascompacts}) we are able to define 
\begin{equation}
 f_{ij}=(\pi_1|_I)^{-1}(e_{ij})\in I\quad \textrm{ for }\quad 1\leq i,j\leq N. \label{eq:matrixunitsinI}
\end{equation}
Set $F=\{ x_1,...,x_s\}\cup \{ f_{i,j}:i,j=1,...,N\}\subset A.$ Choose $\delta>0$ satisfying:
\begin{equation}
 \frac{1}{1+\delta}-316\delta^{1/8}>0\quad \textrm{ and }\quad N\cdot 90 \delta^{1/512}<\epsilon/2.  \label{eq:deltaconditions}
\end{equation}
Now choose $\eta>0$ such that $\eta<\delta$ and $2\eta<\delta(\delta,N)$ (as defined in Theorem \ref{thm:spllitingupwhenrangeisarbitrary}).
 Apply Lemma \ref{cor:usefulUCPOL} to obtain a finite dimensional $C^*$-algebra $B$ and a UCP map $\phi:B\rightarrow A$  such that (without loss of generality)
\begin{equation}
F\subset \phi(B)\quad  \textrm{and}\quad \V\phi^{-1}\V_{cb}<1+\eta<1+\delta, \label{eq:Falmostcontainedinrange} 
\end{equation}
 where  ``without loss of generality'' pertains to taking $F\subset\phi(B)$ instead of $\phi(B)$ close to $F.$

Our first step is to obtain a minimal central projection $p_\ell\in B$ such that $\pi_1\phi|_{p_\ell B}$ can be perturbed to a complete order embedding.  In turn, this will produce a finite rank projection satisfying (\ref{eq:QDconditions}).

Use Wittstock's extension Theorem to obtain a completely bounded extension $\psi:B(H_1)\oplus\cdots\oplus B(H_n)\rightarrow B$  of $\phi^{-1}$ with $\V\psi\V_{cb}\leq 1+\eta.$  Then $\psi$ is unital and  by replacing $\psi$ with $1/2(\psi+\psi^*)$ if necessary, we may assume that $\psi$ is self-adjoint.  
We use  \cite[Proposition 1.19]{Wassermann94} applied to $\psi$ to obtain a UCP map 

\begin{equation}
 T:B(H_1)\oplus\cdots\oplus B(H_n)\rightarrow B\quad \textrm{ with }\quad \V T\pi\phi-id_B\V_{cb}\leq\eta. \label{eq:Tthatinvertsphi}
\end{equation}
By the definition of $I$, we have 
\begin{equation}
 pB(H_1)p\subset \pi_1(\phi(B)\cap I)\oplus 0\oplus \cdots\oplus 0=\pi(\phi(B)\cap I). \label{eq:inclusioninpiI}
\end{equation}
Let $S$ be the restriction of $T$ to $pB(H_1)p.$  It follows from (\ref{eq:Tthatinvertsphi}) and (\ref{eq:inclusioninpiI}) that
\begin{equation}
 S \quad\textrm{ is an injective CPC, with }\quad\V S^{-1}\V_{cb}\leq (1-\eta)^{-1}<1+2\eta \label{eq:conditionsonS}
\end{equation}

We combine (\ref{eq:matrixunitsinI})-(\ref{eq:conditionsonS})  in the following ``approximately'' commutative diagram:
\begin{equation}
\xymatrix@!=2.5pc{
\pi_1(I)  && I &\\
pB(H_1)p \ar@/_3pc/[rrrr]^S_{\sim  2\eta} \ar[u]_\cup \ar[rr]^{(\pi_1|_I)^{-1}} && \phi(B)\cap I \ar[u]_\cup \ar[rr]^-{\phi^{-1}} && B 
 } 
\end{equation}
Since $B$ is finite dimensional, there are integers $n_1,...,n_k$ such that
\begin{equation*}
 B\cong \bigoplus_{i=1}^k M_{n_i}\quad\textrm{ with minimal central projections }\quad p_1,...,p_k\in B.
\end{equation*}
Recall that $2\eta<\delta(\delta,N)$, so by (\ref{eq:conditionsonS}) and  Theorem \ref{thm:spllitingupwhenrangeisarbitrary}, there is an index $1\leq\ell\leq k$ such that $S_{p_\ell}$ is injective with
\begin{equation}
\V S_{p_{\ell}}^{-1}\V_{cb}\leq 1+\delta. \label{eq:cutdownSinjective}
\end{equation}
Recall $f_{11},...,f_{NN}$ as defined in (\ref{eq:matrixunitsinI}). By Corollary \ref{cor:liftingrank1projectionsgeneral},
there exist rank 1 projections $r_1,...,r_N\in B$ such that $\V \phi(r_i)-f_{ii}\V\leq 315\delta^{1/8}$, for $i=1,...,N.$ 
We have
\begin{align}
 \V p_\ell r_i-p_\ell S\pi_1(f_{ii})\V &\leq \V r_i-S\pi_1(f_{ii})\V \notag \\
&\leq\V T\pi(\phi(r_i))-T\pi(f_{ii})\V+\delta \quad \textup{(by (\ref{eq:Tthatinvertsphi}))} \notag \\
&\leq \V \phi(r_i)-f_{ii}\V+\delta\leq 315\delta^{1/8}+\delta. \label{align:anotherperturbation}
\end{align}
By (\ref{align:anotherperturbation}), (\ref{eq:cutdownSinjective}) and (\ref{eq:deltaconditions}) we have for $i=1,...,N$
\begin{equation*}
 \V p_\ell r_i\V\geq \V S_{p_\ell}\pi_1(f_{ii})\V-316\delta^{\frac{1}{8}}\geq \frac{1}{1+\delta}-316\delta^{\frac{1}{8}}>0.
\end{equation*}
Since each $r_i$ is a rank 1 projection, it follows that
\begin{equation}
r_i\in M_{n_\ell}\quad \textrm{ for }i=1,...,N. \label{eq:locationofr_i}
\end{equation}
At this point, the map $S$ is no longer useful to us,  we only needed it to prove (\ref{eq:locationofr_i}).  We will now show that $\pi_1\phi|_{p_\ell B}$ can be perturbed to a complete order embedding.  To this end, first note that

\begin{equation*}
 \V (\pi_2\oplus\cdots\oplus \pi_n)(\phi(r_1))\V\leq \V (\pi_2\oplus\cdots\oplus \pi_n)(f_{11})\V+315\delta^{1/8}=315\delta^{1/8}.
\end{equation*}
Since $\V (\pi\phi)^{-1}\V_{cb}<1+\delta$, by \cite[Corollary 2.7]{Eckhardt08}  it follows that $\pi_1\phi|_{p_\ell B}$ is injective with 
$\V(\pi_1\phi|_{p_\ell B})^{-1}\V_{cb}\leq 1+3\delta^{1/3}.$ Since $e_{11}=\pi_1(f_{11})$ is a rank 1 projection and $\V \pi_1\phi(r_1)-e_{11}\V\leq \V \phi(r_1)-f_{11}\V\leq 315\delta^{\frac{1}{8}}$, we   apply Theorem \ref{thm:minimalperturb} to obtain a complete order embedding $\psi:p_\ell B\cong M_{n_\ell}\rightarrow B(H_1)$ such that
\begin{equation*}
 \V \psi-\pi_1\phi|_{p_\ell B}\V_{cb}\leq 1360(315\delta^{\frac{1}{8}})^{\frac{1}{32}}\leq 2000\delta^{1/256}.
\end{equation*}
Moreover, by Theorem \ref{thm:homosaredense} we may assume there is a rank $n_\ell$ projection $q\in B(H_1)$ such that
$\psi=\psi_q+\psi_{(1-q)}$ with $\psi_q$ a nonzero *-homomorphism.

To finish this portion of the proof we will show that for every $x\in F$ we have
\begin{equation*}
 \V [\pi_1(x), q]\V\leq \epsilon\quad\textrm{ and }\quad \V q\pi_1(x)q\V\geq1-\epsilon.
\end{equation*}
Since $F$ is  contained  in  $\phi(B)$, the first inequality follows from Lemma \ref{lem:usefulcalculations}.(\ref{eq:cutdownalmostcommute}). For the second inequality, let $1\leq i\leq N,$  then,
\begin{align*}
 \V e_{ii}qe_{ii}\V &= \V qe_{ii}q\V\geq \V q\pi\phi(r_i)q\V-315\delta^{1/8}\\
&\geq \V q\psi(r_i)q\V-2001\delta^{1/256}=1-2001\delta^{1/256}.
\end{align*}
Since $e_{ii}$ is a rank 1 projection, it follows that 
\begin{equation}
\V e_{ii}-e_{ii}q\V^2=\V e_{ii}-e_{ii}qe_{ii}\V \leq 2001\delta^{1/256}.\label{eq:e_iiclosetocutdwonq}
\end{equation}
Recall  from (\ref{eq:QDconditions}) that $p=\sum_{i=1}^N e_{ii}.$ By applying (\ref{eq:e_iiclosetocutdwonq}), we have
\begin{equation*}
 \V p-pq\V\leq N\cdot 45 \delta^{1/512}.
\end{equation*}
Finally, by (\ref{eq:QDconditions}), for every $x\in F$  we have
\begin{equation*}
 \V q\pi_1(x)q\V \geq \V pq\pi_1(x)qp\V \geq \V p\pi_1(x)p\V -N\cdot 90 \delta^{1/512}\geq 1-\epsilon/2-N\cdot 90 \delta^{1/512}\geq1-\epsilon.
\end{equation*}
Recall from (\ref{eq:QDconditions}) that $p$ dominates the arbitrary finite rank projection $r.$  By \cite[Proposition 3.6]{Brown04} it follows that $\pi_1$ is a quasidiagonal representation of $A.$ This completes the proof of Lemma \ref{lem:GCRQDrep}.
\newline
We now return to the proof of the Theorem.  All variables defined in the proof of Lemma \ref{lem:GCRQDrep} are now free.

Define $\pi_{GCR}=\pi_1\oplus\cdots\oplus \pi_r.$ Let $F$ be a finite subset of the unit sphere of $A$ and $\epsilon>0.$
If all of the $J_i$ are GCR ideals, then $A$ is inner quasidiagonal by Lemma \ref{lem:GCRQDrep}.  If not, then assume that there is an $x\in F$ such that $\V \pi_{GCR}(x)\V<\V x\V.$ Let
\begin{equation*}
 G:=\{ y\in F: \V y\V>\V \pi_{GCR}(y)\V\}\quad \textrm{ and }\quad \gamma:=\min \{\V y\V-\V \pi_{GCR}(y)\V:y\in G\}>0.
\end{equation*}
Choose $\delta>0$ such that
\begin{equation}
 \delta<\Big(\frac{\epsilon}{|F|^{3/2} 213}\Big)^{12n}\quad\textrm{ and }\quad \frac{1}{1+3\delta^{\frac{2}{3n}}}-3\delta>1-\gamma.                \label{eq:deltaconditionspart2}
\end{equation}

Use Lemma \ref{cor:usefulUCPOL} to obtain a finite dimensional $C^*$-algebra $B$ and a UCP map $\phi:B\rightarrow A$
such that $\V \phi^{-1}\V_{cb}<1+\delta$ and (without loss of generality) $F\cup\{yy^*:y\in F\}\subset \phi(B).$

Since $B$ is finite dimensional, there are integers $n_1,...,n_k$ such that

\begin{equation*}
 B\cong \bigoplus_{i=1}^k M_{n_i}\quad\textrm{ with minimal central projections }\quad p_1,...,p_k\in B.
\end{equation*}
For each $y\in G$ choose an index $1\leq i_y\leq k$ such that $\V p_{i_y}\phi^{-1}(y)\V=\V \phi^{-1}(y)\V$ and let
\begin{equation*}
 M=\bigcup_{y\in G} \{i_y\}\subset \{1,\cdots,k\}.
\end{equation*}
For $i\in M$ and $1\leq j\leq n$ define the maps $\phi_{i,j}:p_iB\rightarrow B(H_j)$ as 
\begin{equation*}
 \phi_{i,j}(x)=\pi_j\phi(x)\quad \textrm{ for all }\quad x\in p_iB.
\end{equation*}

  By repeated use of \cite[Corollary 2.7]{Eckhardt08}, for each $i\in M$ there is an index $1\leq j(i)\leq n$ such that $\phi_{i, j(i)}$ is injective with 
\begin{equation}
 \V \phi_{i, j(i)}^{-1}\V_{cb}\leq 1+3\delta^{\frac{2}{3n}}. \label{eq:boundsoncoordinates}
\end{equation}
We now show that $j(i)>r$ for $i\in M.$ By \cite[Proposition 1.19]{Wassermann94} there is a UCP map $T:A\rightarrow B$ such that $\V T|_{\phi(B)}-\phi^{-1} \V_{cb}\leq\delta.$ Let $i\in M$ and $y\in G$ such that $i=i_y.$  By (\ref{eq:boundsoncoordinates}) we have
\begin{align*}
 \frac{1}{1+3\delta^{\frac{2}{3n}}}&\leq \frac{\V \phi^{-1}(y)\V^2}{1+3\delta^{\frac{2}{3n}}}=\frac{\V p_i\phi^{-1}(y)\phi^{-1}(y)^*\V}{1+3\delta^{\frac{2}{3n}}}\\
&\leq \V \phi_{i, j(i)}( p_i\phi^{-1}(y)\phi^{-1}(y)^*)\V\\
&= \V \pi_{j(i)}\phi(p_i\phi^{-1}(y)\phi^{-1}(y)^*)\V\\
&\leq \V \pi_{j(i)}\phi(\phi^{-1}(y)\phi^{-1}(y)^*)\V\\
&\leq \V \pi_{j(i)}\phi(T(y)T(y)^*)\V+2\delta\\
&\leq \V \pi_{j(i)}\phi(T(yy^*))\V+2\delta\\
&\leq \V \pi_{j(i)}\phi(\phi^{-1}(yy^*))\V+3\delta\\
&\leq \V \pi_{j(i)}(yy^*)\V+3\delta.
\end{align*}
Hence $\V\pi_{j(i)}(yy^*)\V>1-\gamma\geq \V\pi_{GCR}(yy^*)\V$  by (\ref{eq:deltaconditionspart2}). Hence $j(i)>r.$ Therefore $\pi_{j(i)}$ is not a GCR representation so $\pi_{j(i)}(A)\cap K(H_{j(i)})=\{0\}$ for each $i\in M.$ We now apply Theorem \ref{thm:antiliminalperturb} to obtain complete order embeddings $\psi^i:M_{n_i}\rightarrow B(H_{j(i)})$ with 
\begin{equation*}
 \V \phi_{i, j(i)}-\psi^i\V_{cb}\leq 273\delta^{1/6n} \quad \textrm{ for each }i\in M.
\end{equation*}
By Theorem \ref{thm:homosaredense} we may assume that for each $i\in M$ there is a rank $n_i$ projection $q_i\in B(H_{j(i)})$ such that 
\begin{equation}
 \psi^i=\psi^i_{q_i}+\psi^i_{(1-q_i)}\quad \textrm{ with }\quad \psi^i_{q_i}\quad\textrm{ a nonzero *-homomorphism.}
\end{equation}
Let $y\in G$ and let $i_y=i\in M.$  Then
\begin{align}
 \V q_i \pi_{j(i)}(y)q_i\V&= \V q_i \pi_{j(i)}(\phi(p_i\phi^{-1}(y))q_i+q_i\pi(\phi((1-p_i)\phi^{-1}(y))q_i\V \notag\\
&\geq \V q_i \phi_{i, j(i)}(p_i\phi^{-1}(y))q_i\V-8(273)\delta^{1/6n} \quad (\textrm{ by Lemma }\textup{\ref{lem:usefulcalculations}.(\ref{eq:Lshapednormcontrol})})\notag\\
&\geq \V q_i \psi^i(p_i\phi^{-1}(y))q_i\V-9(273)\delta^{1/6n}\notag\\
&=\V p_i\phi^{-1}(y)\V-9(273)\delta^{1/6n}\notag\\
&\geq 1-9(273)\delta^{1/6n}\geq 1-50\delta^{1/12n}.\label{align:normcutdownbyq_i}
\end{align}
Furthermore, for any $x\in F$ we have by Lemma \ref{lem:usefulcalculations}.(\ref{eq:cutdownalmostcommute}), that
\begin{equation}
 \V [\pi(x), q_i]\V \leq 8(273\delta^{1/6n})^{1/2}\leq 133\delta^{1/12n}. \label{eq:Fapproxcommuteswithq_i}
\end{equation}
Now we show that the sum of the $q_i$ is almost a projection.
\begin{align*}
 \Big\V \sum_{i\in M}q_i\Big\V&=\Big\V \sum_{i\in M}q_i\psi^i(p_i)q_i\Big\V\leq \Big\V \sum_{i\in M}\psi^i(p_i)\Big\V\leq \Big\V \sum_{i\in M}\phi_{i, j(i)}(p_i)\Big\V+|M|273\delta^{1/6n}\\
&\leq \Big\V \sum_{i\in M}\phi(p_i)\Big\V+|M|273\delta^{1/6n}\leq 1+|M|273\delta^{1/6n}.
\end{align*}
Then for each $i\in M$ we have
\begin{align*}
 \Big\V q_i\Big(\sum_{j\neq i}q_j\Big)\Big\V&\leq \Big\V q_i\Big(\sum_{j\neq i}q_j\Big)^{1/2} \Big\V(1+|M|273\delta^{1/6n})^{1/2}\\
&=\Big\V q_i\Big(\sum_{j\neq i}q_j\Big)q_i \Big\V^{1/2}(1+|M|273\delta^{1/6n})^{1/2}\\
&\leq (|M|273\delta^{1/6n})^{1/2}(1+|M|273\delta^{1/6n})^{1/2}\leq 40|M|^{1/2}\delta^{1/12n}.
\end{align*}
From which it follows that
\begin{equation*}
 \Big\V \sum_{i\in M} q_i-\Big(\sum_{i\in M} q_i\Big)^2\Big\V\leq 40|M|^{3/2}\delta^{1/12n}.
\end{equation*}

Hence, basic spectral theory produces a projection $q\in B(H_{r+1})\oplus\cdots\oplus B(H_n)$ such that
\begin{equation}
 \Big\V q-\sum_{i\in M} q_i\Big\V\leq 80|M|^{3/2}\delta^{1/12n}.  \label{eq:projectionclosetoq_isum}
\end{equation}
Combining (\ref{align:normcutdownbyq_i}) and (\ref{eq:projectionclosetoq_isum}), it follows that
\begin{equation}
 \V q \pi_{j(i)}(y)q\V\geq 1-|M|^{3/2} 130\delta^{1/12n}\geq 1-\epsilon \label{eq:qapproxpreservenorm}
\end{equation}
and then by (\ref{eq:Fapproxcommuteswithq_i}) and (\ref{eq:projectionclosetoq_isum}) it follows that
\begin{equation}
 \V [\pi(x), q]\V \leq 213|M|^{3/2}\delta^{1/12n}\leq\epsilon \quad\textrm{ for all }\quad x\in F.\label{eq:qapproxcommutes}
\end{equation}
We use Lemma \ref{lem:GCRQDrep} to obtain a finite rank projection $p\in B(H_1)\oplus\cdots\oplus B(H_r)$ such that
\begin{equation}
 \max_{x\in F}\V [p, \pi(x)]\V \leq \epsilon\quad \textrm{ and }\min_{x\in F\setminus G}\V p\pi(x)p\V\geq 1-\epsilon.  \label{eq:projectionfromGCR}
\end{equation}
Finally combining (\ref{eq:qapproxcommutes}) , (\ref{eq:qapproxpreservenorm}) and (\ref{eq:projectionfromGCR}) if follows that
\begin{equation*}
\max_{x\in F}\V [p+q, \pi(x)]\V \leq \epsilon\quad \textrm{ and }\min_{x\in F}\V (p+q)\pi(x)(p+q)\V\geq 1-\epsilon.
\end{equation*}
Since $p$ and $q$ are orthogonal, it follows that  $A$ is inner quasidiagonal.

\end{proof}
\subsection{Behavior of $\mathcal{OL}_\infty$}
Theorem \ref{thm:maintheoremin3} also provides  somewhat surprising results about the behavior of the invariant $\mathcal{OL}_\infty$ that we now describe.  Recall from \cite[Proposition 4.4]{Eckhardt08}, that for any nuclear $C^*$-algebras $A$ and $B$ we have
\begin{equation}
 \mathcal{OL}_\infty(A\otimes B)\leq \mathcal{OL}_\infty(A)\mathcal{OL}_\infty(B).\label{eq:tensorproductmult}
\end{equation}
This inequality can be strict in general.  Consider the Cuntz algebra $\mathcal{O}_2.$  By \cite[Theorem 3.4]{Junge03}, it follows that $\mathcal{OL}_\infty(\mathcal{O}_2)>\sqrt{(1+\sqrt{5})/2}>1.$  By a result of George Elliott (see \cite[Chapter 5]{Rordam02}) we have $\mathcal{O}_2\otimes \mathcal{O}_2\cong \mathcal{O}_2,$ showing that (\ref{eq:tensorproductmult}) can be strict.

It is at least reasonable to think that equality in (\ref{eq:tensorproductmult}) will hold if one of the algebras is strong NF. We now show that (\ref{eq:tensorproductmult}) can be strict even when one of the algebras is AF.

 As pointed out in \cite{Blackadar01}, there is a nuclear, primitive quasidiagonal $C^*$-algebra $\mathcal{B}$ (defined in  \cite{Brown84}) that is not inner quasidiagonal, and hence not a strong NF algebra.  By Theorem \ref{thm:maintheoremin3}, it follows that $\mathcal{OL}_\infty(\mathcal{B})>1.$  Let $U$ be any UHF algebra. Then $U$ is a strong NF algebra and $\mathcal{B}\otimes U$ is primitive, quasidiagonal and antiliminal.  By \cite[Corollary 2.6]{Blackadar01} it follows that $\mathcal{B}\otimes U$ is a strong NF algebra.  Therefore,
\begin{equation*}
 \mathcal{OL}_\infty(\mathcal{B}\otimes U)=1<\mathcal{OL}_\infty(\mathcal{B}).
\end{equation*}
The fact that $U$ is antiliminal played a key role in this example.  Therefore, it would be interesting to know what happens for the Type I building blocks.  In particular,
\begin{question} Let $B$ be the compact operators or a commutative $C^*$-algebra and $A$ a nuclear $C^*$-algebra.  Do we have
\begin{equation*}
 \mathcal{OL}_\infty(A\otimes B)=\mathcal{OL}_\infty(B)?
\end{equation*}
\end{question}
We now show, using the same example, that $\mathcal{OL}_\infty$ is not continuous with respect to inductive limits.  We need the following result:
\begin{theorem} \label{thm:OLpAp} Let $A$ be a unital $C^*$-algebra with $\mathcal{OL}_\infty(A)=1$  and let $p\in A$ be a projection.  Then $\mathcal{OL}_\infty(pAp)=1.$
\end{theorem}
\begin{proof} Let $p\in F\subset pAp$ be a finite subset and $\delta>0$.
Then obtain a finite dimensional $C^*$-algebra $B$ and  maps $\phi:B\rightarrow A$ and $\psi:A\rightarrow B$  that satisfy Lemma \ref{cor:usefulUCPOL}
for $F$ and $\delta>0.$ Since $\psi$ is $\delta$-multiplicative on $F$, it follows that there is a projection $q\in B$ such that $\V\psi(p)-q\V\leq 2\delta.$  Then for every $x\in (qBq)^+$ of norm 1, we have $\phi(x)\leq \phi(q).$
From this it follows that for every $x\in qBq$ we have
\begin{equation*}
 \V p\phi(x)p-\phi(x)\V\leq 12\delta.
\end{equation*}
Finally, we apply \cite[Lemma 2.13.2]{Pisier03} to obtain a map $\alpha:qBq\rightarrow pAp$  with
\begin{equation*}
 \V \alpha\V_{cb}\V \alpha^{-1}\V_{cb}\leq \frac{1+12\delta}{1-12\delta}\quad \textrm{ and } F \subset \alpha(qBq).\end{equation*} 
\end{proof}
\begin{corollary} \label{cor:uniformmatrixbound} Let $A$ be a unital $C^*$-algebra with $\mathcal{OL}_\infty(A)>1.$  Then there exists a constant $r>1$ such that
\begin{equation}
 \inf_{n\in\mathbb{N}} \mathcal{OL}_\infty(M_n\otimes A)\geq r.
\end{equation}
\end{corollary}
\begin{proof}
If not, then we could apply the proof of Theorem  \ref{thm:OLpAp} to the sequence of algebras $e_{11}\otimes A\subseteq M_n\otimes A$ and deduce that $\mathcal{OL}_\infty(A)=1.$
\end{proof}
It is obvious from the definition of $\mathcal{OL}_\infty$ that if $A$ is the inductive limit of the nuclear $C^*$-algebras $A_n$, then
\begin{equation*}
 \limsup_{n\rightarrow\infty} \mathcal{OL}_\infty(A_n)\geq \mathcal{OL}_\infty (A).
\end{equation*}
We show that this inequality can be strict.  Consider the algebras $\mathcal{B}$  and  $U$ as above.   Since $U$ is the inductive limit of matrix algebras $(M_{n_k})_{k=1}^\infty$, it follows that $U\otimes \mathcal{B}$ is the inductive limit of the algebras $M_{n_k}\otimes\mathcal{B}.$ By Corollary \ref{cor:uniformmatrixbound},
\begin{equation*}
 \limsup_{k\rightarrow\infty} \mathcal{OL}_\infty(M_{n_k}\otimes\mathcal{B})>1=\mathcal{OL}_\infty (U\otimes \mathcal{B}).
\end{equation*}
\section{Concluding Remarks}
Recall from the Introduction the following question:
\begin{question} Let A be a $C^*$-algebra with $\mathcal{OL}_\infty(A)=1.$  Is A a strong NF algebra? 
\end{question}
This question is still open, but we note the following necessary conditions:
\begin{theorem} \label{thm:finalremarksoni}
 Suppose A is a $C^*$-algebra with $\mathcal{OL}_\infty(A)=1,$ but A is not a strong NF algebra.  Then
\begin{enumerate}
 \item A does not have a separating family of irreducible quasidiagonal representations.
\item A is quasidiagonal.
\item A has a separating family of irreducible stably finite representations.
\item For every finite subset $F\subseteq \textup{Prim}(A)$, $\cap_{J\in F}J\neq\{0\}.$
\end{enumerate}
\end{theorem}
 \begin{proof}
The conclusions follow from \cite[Corollary 1.3]{Blackadar07}, \cite[Theorem 3.2]{Junge03},  \cite[Theorem 5.4]{Eckhardt08} and Theorem \ref{thm:maintheoremin3} respectively.
 \end{proof}
We don't have to look very far to find a nuclear $C^*$-algebra that satisfies 1-4 in Theorem \ref{thm:finalremarksoni}.
 Recall the $C^*$-algebra $\mathcal{B}$ from above.  Then  $C[0,1]\otimes \mathcal{B}$  satisfies 1-4 in Theorem \ref{thm:finalremarksoni}.  Therefore, it would certainly be interesting to know $\mathcal{OL}_\infty(C[0,1]\otimes\mathcal{B}).$  Unfortunately, we have been unable to decide this question.
\section*{Acknowledgment} This work constitutes a portion of the author's Ph.D. thesis at the University of Illinois at Urbana-Champaign.  I would like to acknowledge the support of the mathematics department and especially my advisor Zhong-Jin Ruan.

\bibliographystyle{plain}
\bibliography{mybib}

\end{document}